\title{Hyperbolic manifolds and tessellations of type \(\{3,5,3\}\) associated with \(L_2(q)\)}

\author{Gareth A. Jones\\
School of Mathematics\\
University of Southampton\\
Southampton SO17  1BJ, U.K.\\
{\tt G.A.Jones@maths.soton.ac.uk}
\and  Cormac D.~Long\\
School of Mathematics\\
University of Southampton\\
Southampton SO17  1BJ, U.K.\\
{\tt clong@soton.ac.uk}
\and   Alexander D.~Mednykh \\
Sobolev Institute for Mathematics\\
Novosibirsk State University\\
Novosibirsk, Russia\\
{\tt smedn@mail.ru}}

\documentclass[12pt]{article}
\usepackage[latin1]{inputenc}
\usepackage{a4wide}
\usepackage{latexsym}
\usepackage{amsmath}
\usepackage{amssymb}
\usepackage{amsfonts}
\newtheorem{thm}{Theorem}[section]

\newtheorem{cor}[thm]{Corollary}

\date{}

\usepackage{colortbl}
\usepackage{rotating}

\begin{document}

\maketitle

\begin{abstract}
\noindent We classify the normal subgroups $K$ of the tetrahedral group $\Delta=[3,5,3]^+$, the even subgroup of the Coxeter group $\Gamma=[3,5,3]$, with $\Delta/K$ isomorphic to a finite simple group $L_2(q)$. We determine their normalisers $N(K)$ in the isometry group of hyperbolic $3$-space ${\cal H}^3$, the isometry groups $N(K)/K$ of the associated hyperbolic $3$-manifolds ${\cal H}^3/K$, and the symmetry groups $N_{\Gamma}(K)/K$ of the icosahedral tessellations of these manifolds, giving a detailed analysis of how $L_2(q)$ acts on these tessellations.
\end{abstract}

%% 20F55 = Reflection and Coxeter groups,
%% 22E40 = Discrete subgroups of Lie groups,
%% 57M60 = Group actions in low dimensions


{\bf Running head:} Hyperbolic manifolds\\

\section{Introduction}

Gehring, Marshall and Martin, building on a series of earlier results in~\cite{GMM98, GMM02, GM98, GM99, GM08}, have recently announced a proof of the long-standing conjecture that among the discrete groups of isometries of hyperbolic $3$-space ${\cal H}^3$, the group of least covolume is the normaliser $\Omega$ of the Coxeter group $\Gamma=[3,5,3]$. This conjecture, proved earlier for arithmetic isometry groups by Chinburg and Friedman in~\cite{CF}, has directed attention towards the torsion-free normal subgroups $K$ of finite index in $\Omega$ and in its orientation-preserving subgroup $\Omega^+$, since these subgroups uniformise the compact hyperbolic $3$-manifolds ${\cal H}^3/K$ which have large isometry groups $\Omega/K$ or orientation-preserving isometry groups $\Omega^+/K$, in the sense of maximising the number of isometries per unit volume. These isometry groups are $3$-dimensional analogues of the Hurwitz groups which arise in the $2$-dimensional case as groups attaining Hurwitz's upper bound of $84(g-1)$ for the number of automorphisms of a compact Riemann surface of genus $g\geq 2$. The finite quotients of $\Omega$ and of $\Omega^+$ are closely related to those of the orientation-preserving subgroup $\Delta=\Gamma^+$ of $\Gamma$, and since $\Delta$ is perfect its nontrivial quotients are all coverings of nonabelian simple groups. By far the most tractable and ubiquitous of the nonabelian finite simple groups are those of the form $L_2(q)=SL_2(q)/\{\pm I\}$ for prime powers $q=p^n\geq 4$. Macbeath~\cite{Mac} has determined which of the groups $L_2(q)$ are Hurwitz groups, and Singerman~\cite{Sin} has extended this to $H^*$-groups, similarly associated with non-orientable Klein surfaces.

Motivated by these results, several authors have recently considered the simple groups $L_2(q)$ and other closely related groups as quotients in the $3$-dimensional context. Paoluzzi~\cite[Theorem~1]{Pao} has determined which groups $L_2(q)$ are quotients of $\Delta$ (her tetrahedral group $T_{5,2}$), and has also obtained partial results on quotients of $\Omega^+$ and $\Gamma$ (her $\tilde T$ and $C_{5,2}$). Torstensson~\cite{Tor} has extended these results by determining all the quotients of $\Omega^+$ (her group $\Gamma$) of type $L_2(q)$ or $PGL_2(q)$. Conder, Martin and Torstensson~\cite{CMT} have also considered the groups $L_2(q)$, and also alternating and symmetric groups, as quotients. The emphasis in these papers has been on the purely algebraic problem of determining those $q$ for which groups such as $L_2(q)$ or $PGL_2(q)$ arise  as quotients of $\Omega, \Omega^+$ or $\Gamma$. However, from the geometric point of view there is also interest in going further to determine all the normal subgroups $K$ with such a quotient, their conjugacy under isometries, since this corresponds to isometry of the corresponding quotient manifolds, and their normalisers $N(K)$ in ${\rm Iso}\,{\cal H}^3$, since these determine the isometry groups $N(K)/K$ of these manifolds. Similarly, results about normality, conjugacy and quotients in $\Gamma$ give combinatorial information about the tessellations ${\cal I}/K$ of these manifolds inherited from the icosahedral tessellation ${\cal I}=\{3,5,3\}$ of ${\cal H}^3$ associated with $\Gamma$.

Jones and Mednykh~\cite{JM} have shown that there is a unique torsion-free normal subgroup $K_0$ of least index in $\Omega$; this has quotient $\Omega/K_0\cong PGL_2(11)\times C_2$, the isometry group of a $3$-manifold ${\cal H}^3/K_0$ tessellated by eleven icosahedra to form the honeycomb $\{3,5,3\}_6$ described by Coxeter in~\cite{Cox70}. By contrast, there are two torsion-free normal subgroups $K_1$ and $K_2$ of least index in $\Omega^+$; these subgroups, which are conjugate in $\Omega$, have quotient $\Omega^+/K_i\cong PGL_2(3^2)$, the isometry group of a chiral pair of $3$-manifolds ${\cal H}^3/K_i$ tessellated by six icosahedra to form the twisted honeycombs $\{3,5,3\}_4$ and $\{3,5,3\}_5$. One can regard these manifolds ${\cal H}^3/K_i\;(i=0, 1, 2)$ as $3$-dimensional analogues of Klein's quartic curve, the Riemann surface of genus $3$ corresponding to the smallest Hurwitz group $L_2(7)$. In this paper we determine all the normal subgroups (necessarily torsion-free) $K$ of $\Delta$ with $\Delta/K\cong L_2(q)$ for any $q$, together with their normalisers $N(K)$, and we describe the resulting isometry groups $N(K)/K$ and symmetry groups $N_{\Gamma}(K)/K$. Some of the methods used here are similar to those applied in~\cite{JL} to the Coxeter group $[5,3,5]$ and its normaliser, associated with the dodecahedral tessellation $\{5,3,5\}$ of ${\cal H}^3$.

Our main result is the following, where $F_q$ denotes the finite field of order $q$:

% Thm 1.1.

\begin{thm} The group $\Delta=[3,5,3]^+$ has only the following normal subgroups $K$ with $\Delta/K\cong L_2(q)$, where $q$ is a power of a prime $p$:
\vskip2pt
\noindent{\rm(a)} if $p=2$ there is one normal subgroup $K$ with $\Delta/K\cong L_2(2^4)$;
\vskip2pt
\noindent{\rm(b)} if $p=5$ there is one normal subgroup $K$ with $\Delta/K\cong L_2(5^2)$;
\vskip2pt
\noindent{\rm(c)} if $p=11$ there is one normal subgroup $K=K_0$ with $\Delta/K\cong L_2(11)$ and one normal subgroup $K$ with $\Delta/K\cong L_2(11^2)$;
\vskip2pt
\noindent{\rm(d)} if $p\equiv\pm 1$ {\rm mod}~$(5)$ with $p\equiv 1, 3, 4, 5$ or $9$ {\rm mod}~$(11)$ there are four normal subgroups $K$ with $\Delta/K\cong L_2(p)$ or two with $\Delta/K\cong L_2(p^2)$ as $3\pm 2\sqrt 5$ are both squares or both non-squares in $F_p$;
\vskip2pt
\noindent{\rm(e)} if $p\equiv\pm 1$ {\rm mod}~$(5)$ with $p\equiv 2, 6, 7, 8$ or $10$ {\rm mod}~$(11)$ there are two normal subgroups $K$ with $\Delta/K\cong L_2(p)$ and one with $\Delta/K\cong L_2(p^2)$;
\vskip2pt
\noindent{\rm(f)} if $2\neq p\equiv\pm 2$ {\rm mod}~$(5)$ and $3\pm 2\sqrt 5$ are both squares in $F_{p^2}$ there are two normal subgroups $K$ with $\Delta/K\cong L_2(p^2)$;
\vskip2pt
\noindent{\rm(g)} if $2\neq p\equiv\pm 2$ {\rm mod}~$(5)$ and $3\pm 2\sqrt 5$ are both non-squares in $F_{p^2}$ there is one normal subgroup $K$ with $\Delta/K\cong L_2(p^4)$.
\end{thm}

The smallest of these quotients is $L_2(3^2)\cong A_6$, of order $360$, corresponding to two normal subgroups $K=K_1$ and $K_2$ in case~(f), and the second smallest is $L_2(11)$, of order $660$, with one normal subgroup $K=K_0$ in case~(c). In fact these are the smallest among {\sl all\/} proper quotients of $\Delta$; this is proved in~\cite{JM}, where the corresponding manifolds and tessellations are studied in some detail.

Note that if $p\equiv\pm 1$ mod~$(5)$ then $5$ is a quadratic residue mod~$(p)$, so $3\pm 2\sqrt 5$ is well-defined as a pair of elements of $F_p$ in case~(d). The significance there of the condition $p\equiv 1, 3, 4, 5$ or $9$ mod~$(11)$ is that these are the odd primes for which $-11$ is a quadratic residue mod~$(p)$; case~(c) also indicates the exceptional status of the prime $11$ here. D\v zambi\'c~\cite{Dza} has explained this, and also the significance of $3\pm2\sqrt 5$ in Theorem~1.1, by regarding $\Delta$ as an arithmetic Kleinian group, associated with a quaternion algebra over the field ${\bf Q}\bigl(\sqrt{3+2\sqrt 5}\bigr)$, and obtaining these subgroups $K$ as the principal congruence subgroups modulo prime ideals dividing $p$; the various cases in Theorem~1.1 correspond to the different prime ideal decompositions of $p$, with $5$ and $11$ (the primes dividing the discriminant $-275$ of this field) the only ramified primes. The special status of the prime $2$ here and in later results is related to the exceptional structure of the groups $L_2(q)$ when $q$ is a power of $2$.

For each of the subgroups $K$ in Theorem~1.1, the group $\Delta/K$ induces a group of orientation-preserving isometries of the manifold ${\cal H}^3/K$, isomorphic to $L_2(q)$. The following result shows that in every case this is a subgroup of index $2$ in a larger group of orientation-preserving isometries induced by $\Omega^+/K$:

% Thm 1.2.

\begin{thm}
If $K$ is a normal subgroup of $\Delta$ with $\Delta/K\cong L_2(q)$ for some $q$ then $K$ is normal in $\Omega^+$ with $\Omega^+/K$ isomorphic to $PGL_2(q)$ or $L_2(q)\times C_2$.
\end{thm}

In \S 4 we give a criterion to determine which of these quotients arises. This is more general than that given by Paoluzzi in~\cite{Pao}, where only the case $q\equiv 1$ mod~$(10)$ is considered, and is rather simpler than that given by Torstensson in~\cite{Tor}. Theorem~1.2 does not extend to all finite simple quotients of $\Delta$: for instance, $\Delta$ has four normal subgroups with quotient isomorphic to the alternating group $A_{25}$, and only two of these are normal in $\Omega^+$.

The Coxeter group $\Gamma=[3,5,3]$ is the symmetry group of a tessellation ${\cal I}=\{3,5,3\}$ of ${\cal H}^3$ by icosahedra with dihedral angles $2\pi/3$. If $K\leq\Gamma$ then ${\cal H}^3/K$ inherits a tessellation ${\cal I}/K$ with symmetry group $N_{\Gamma}(K)/K$. The subgroups $K$ classified in Theorem~1.1 have $N_{\Gamma}(K)=\Gamma$ or $\Delta$ as they are normal in $\Gamma$ or occur in conjugate pairs, and the corresponding tessellations are respectively either reflexible or occur in chiral pairs. It is therefore of interest to decide when such a subgroup $K$ is normal in $\Gamma$, and if so to determine the structure of $\Gamma/K$. In order to state our main result on this topic, we need to introduce some notation. Recall that $P\Sigma L_2(q)$ is the extension of $L_2(q)$ by its group of automorphisms induced by the Galois group ${\rm Gal}\,F_q/F_p\cong C_n$, where $q=p^n$. If $n=4$ we define $P\Sigma L_2(q)^+$ to be the unique subgroup of index $2$ in $P\Sigma L_2(q)$, that is, the extension of $L_2(q)$ by the group of automorphisms induced by ${\rm Gal}\,F_q/F_{p^2}\cong C_2$.

% Thm 1.3.

\begin{thm}
Let $K$ be a normal subgroup of $\Delta$ with $\Delta/K\cong L_2(q)$ for some $q=p^n$;
\vskip2pt
\noindent{\rm(a)} if $p=2$ then $K$ is normal in $\Gamma$ with $\Gamma/K\cong P\Sigma L_2(2^4)^+$.
\vskip2pt
\noindent{\rm(b)} if $q=5^2$ then $K$ is normal in $\Gamma$ with $\Gamma/K\cong P\Sigma L_2(5^2)$;
\vskip2pt
\noindent{\rm(c)} if $q=11$ then $K=K_0$ is normal in $\Gamma$ with $\Gamma/K\cong L_2(11)\times C_2$, and if $q=11^2$ then $K$ is normal in $\Gamma$ with $\Gamma/K\cong P\Sigma L_2(11^2)$;
\vskip2pt
\noindent{\rm(d)} if $p\equiv\pm 1$ {\rm mod}~$(5)$ with $p\equiv 1, 3, 4, 5$ or $9$ {\rm mod}~$(11)$, and $3\pm 2\sqrt 5$ are both squares in $F_p$, then the four normal subgroups $K$ with $\Delta/K\cong L_2(p)$ form two conjugate pairs in $\Gamma$, whereas if $3\pm 2\sqrt 5$ are both non-squares then the two subgroups $K$ with $\Delta/K\cong L_2(p^2)$ are normal in $\Gamma$ with $\Gamma/K\cong P\Sigma L_2(p^2)$;
\vskip2pt
\noindent{\rm(e)} if $p\equiv\pm 1$ {\rm mod}~$(5)$ with $p\equiv 2, 6, 7, 8$ or $10$ {\rm mod}~$(11)$, then the two normal subgroups $K$ with $\Delta/K\cong L_2(p)$ are conjugate in $\Gamma$, whereas the one with $\Delta/K\cong L_2(p^2)$ is normal in $\Gamma$ with $\Gamma/K\cong P\Sigma L_2(p^2)$;
\vskip2pt
\noindent{\rm(f)} if $2\neq p\equiv\pm 2$ {\rm mod}~$(5)$ and $3\pm 2\sqrt 5$ are both squares in $F_{p^2}$, then the two normal subgroups $K$ with $\Delta/K\cong L_2(p^2)$ are conjugate in $\Gamma$;
\vskip2pt
\noindent{\rm(g)} if $2\neq p\equiv\pm 2$ {\rm mod}~$(5)$ and $3\pm 2\sqrt 5$ are both non-squares in $F_{p^2}$, then the normal subgroup $K$ with $\Delta/K\cong L_2(p^4)$ is normal in $\Gamma$ with $\Gamma/K\cong P\Sigma L_2(p^4)^+$.
\end{thm}

The isometry group of the manifold ${\cal H}^3/K$ can be identified with $N(K)/K$. Results of Derevnin and Mednykh~\cite{DM} show that $N(K)\leq\Omega$. Since $\Omega$ is generated by $\Gamma$ and $\Omega^+$, Theorem~1.2 therefore implies the following:

% Coroll 1.4.

\begin{cor}
Let $K$ be a normal subgroup of $\Delta$ with $\Delta/K\cong L_2(q)$ for some $q$. Then the normaliser $N(K)$ of $K$ in ${\rm Iso}\,{\cal H}^3$ is either $\Omega$ or $\Omega^+$ as $K$ is or is not normal in $\Gamma$.
\end{cor}

If $K$ is not normal in $\Gamma$ then Theorem~1.2 and the criterion in \S 4 describe the isometry group $\Omega^+/K$ of ${\cal H}^3/K$. If $K$ is normal in $\Gamma$ then the isometry group is $\Omega/K$, and one can deduce its structure from this information about $\Omega^+/K$ and that about $\Gamma/K$ in Theorem~1.3: if $\Omega^+/K\cong L_2(q)\times C_2$ then $\Omega/K\cong\Gamma/K\times C_2$, whereas if $\Omega^+/K\cong PGL_2(q)$ then $\Omega/K$ is isomorphic to the subgroup of ${\rm Aut}\,L_2(q)=P\Gamma L_2(q)$ generated by $PGL_2(q)$ and the image of $\Gamma/K$.

By Theorem~1.3(c) there is a normal subgroup $N_0$ of $\Gamma$ such that $N_0/K_0$ is the direct factor of $\Gamma/K_0\cong L_2(11)\times C_2$ isomorphic to $C_2$. This subgroup $N_0$ is unique in the following sense:

% Coroll 1.5.

\begin{cor}
The only proper normal subgroups of $\Gamma$ with quotient isomorphic to a subgroup of $L_2(q)$ for any prime power $q$ are $\Delta$ and $N_0$, with quotients $C_2$ and $L_2(11)$.
\end{cor}

Paoluzzi has shown in~\cite{Pao} that $L_2(q)\times C_2$ is a quotient of $\Gamma$ (denoted there by $C_{5,2}$) if and only if $q=11$; Corollary~E strengthens this slightly by implying that $K_0$ is the only normal subgroup of $\Gamma$ with such a quotient. The existence of $N_0$ realises $L_2(11)$ as the symmetry group of Coxeter's $11$-cell (or hendecachoron) ${\cal I}/N_0$~\cite{Cox}. This tessellated orbifold can be regarded as a regular abstract polytope of rank $4$ whose facets are eleven hemi-icosahedra (quotients of icosahedra under antipodal identification). Leemans and Schulte~\cite{LS} have shown that the only other regular abstract polytope of rank greater than $3$ with symmetry group $L_2(q)$ is the $57$-cell (or pentacontaheptachoron), with $57$ hemidodecahedral facets. This has symmetry group $L_2(19)$, arising as a quotient of the Coxeter group $[5,3,5]$; quotients isomorphic to $L_2(q)$ of this and related groups have been studied more generally by Gradolato and Zimmermann~\cite{GZ}, by Paoluzzi~\cite{Pao}, and by two of the present authors~\cite{JL}.

Monson and Schulte~\cite{MS} have used modular reduction of various rank $4$ Coxeter groups $\Gamma$, including $[3,5,3]$ and $[5,3,5]$, to construct polytopes with $4$-dimensional finite orthogonal groups as symmetry groups. In fact, the resulting quotients of $\Gamma^+$ are directly related to quotient groups of type $L_2$: specifically, for odd $q$ one has $O_4^+(q)\cong L_2(q)\times L_2(q)$ and $O_4^-(q)\cong L_2(q^2)$, where the superscript $+$ or $-$ indicates that the corresponding quadratic form has Witt index $2$ or $1$, with similar results for $q=2^e$ (see~\cite[Ch.~II \S 9, \S 10]{Die}, for instance). In particular, when $\Gamma=[3,5,3]$ the corresponding kernels are the normal subgroups $K$ appearing in Theorem~1.1, or intersections of them for a given $q$; a similar remark applies to the discussion of $\Gamma=[5,3,5]$ in~\cite{JL}. The advantage of our more group-theoretic approach is that it allows one to show that all the normal subgroups with a given quotient of type $L_2$ have been obtained.

This paper is organised as follows. In \S 2 we summarise some basic properties of the groups $\Gamma$, $\Delta$ and $\Omega$. These are used in \S\S 3--7, where Theorems~1.1, 1.2 and 1.3 and Corollaries~1.4 and 1.5 are proved. In \S 8 we give a detailed analysis of the action of the various quotients $L_2(q)$ of $\Delta$ on the associated tessellations ${\cal I}/K$, determining the cycle structures of each element on vertices, edges, faces and cells, and also illustrating in the case $q\equiv 11$ mod~$(60)$ how the permutation characters of $L_2(q)$ on these objects can be expressed as sums of irreducible characters. In \S 9 we study the action of $L_2(q)$ on the manifolds ${\cal M}={\cal H}^3/K$, and in \S 10 we show that in the cases $q=11$, $29$ and $59$ (but no others) there is a subgroup $S\leq L_2(q)$ acting regularly on the icosahedra tessellating $\cal M$, so that ${\cal M}/S$ is a manifold consisting of a single icosahedron with its faces identified in pairs. Everitt~\cite{Eve} has classified those manifolds obtained from platonic solids by identification of faces, and we are able to identify each manifold ${\cal M}/S$ with the appropriate entry in his list. In the cases $q=11$ and $29$ we also realise these manifolds ${\cal M}/S$ as cyclic coverings of $S^3$ branched over certain knots.

In \S 11 we give tables of our results for all primes $p\leq 251$, extending those in~\cite{Pao} and~\cite{Tor}; these were obtained by using GAP to check the relevant conditions in finite fields.

\section{Properties of  $\Gamma$ and related groups}

Let $T$ be a tetrahedron in ${\cal H}^3$ with vertices $A, B, C, D$ and dihedral angles $\pi/2$ along the edges $AC, BC$ and $BD$, angles $\pi/3$ along $AB$ and $CD$, and angle $\pi/5$ along $AD$. The Coxeter group $\Gamma=[3,5,3]$ has a presentation
\[\langle a, b, c, d \mid a^2=b^2=c^2=d^2=(ab)^3=(bc)^5=(cd)^3=(ac)^2=(ad)^2=(bd)^2=1\rangle,\]
where $a, b, c$ and $d$ are the reflections of ${\cal H}^3$ in the faces of $T$ opposite $A, B, C$ and $D$. The images of $T$ under $\Gamma$ tessellate ${\cal H}^3$. The images of $T$ under the Coxeter subgroup
\[\Gamma_0=\langle a, b, c\rangle=[3,5]\cong A_5\times C_2\]
of $\Gamma$ form an icosahedron $I$ with dihedral angles $2\pi/3$ and symmetry group $\Gamma_0$, and the images of $I$ under $\Gamma$ form the tessellation ${\cal I}=\{3,5,3\}$ of ${\cal H}^3$, with symmetry group $\Gamma$.

The orientation-preserving subgroup $\Delta=\Gamma^+$ of $\Gamma$ has a presentation
\[
\langle \alpha, \beta, \gamma \mid \alpha^3=\beta^5=\gamma^3=
(\alpha\beta)^2=(\beta\gamma)^2=(\alpha\beta\gamma)^2=1\rangle,
\eqno(2.1)\]
where $\alpha=ab$, $\beta=bc$ and $\gamma=cd$ are rotations of ${\cal H}^3$ through $2\pi/3$, $2\pi/5$ and $2\pi/3$ around the edges $CD, AD$ and $AB$ of $T$. As shown in~\cite{JM}, $\Delta$ is perfect, and every proper normal subgroup of $\Delta$ is torsion-free.
The orientation-preserving subgroup $\Gamma_0^+=\Gamma_0\cap \Delta$ of $\Gamma_0$ is
\[\Delta_0=\langle \alpha, \beta \mid \alpha^3=\beta^5=(\alpha\beta)^2=1\rangle=[3,5]^+\cong A_5.\]

The normaliser $\Omega$ of $\Gamma$ in ${\rm Iso}\,{\cal H}^3$ is a semidirect product of $\Gamma$ by a group $\langle r\rangle\cong C_2$, where $r$ is a half-turn of ${\cal H}^3$ about an axis through the mid-points of the edges $AD$ and $BC$ of $T$; its action by conjugation on $\Gamma$ is given by transposing $a$ and $d$, and transposing $b$ and $c$. The orientation-preserving subgroup $\Omega^+$ of $\Omega$ is a semidirect product of $\Delta$ by $\langle r\rangle$, with $r$ transposing $\alpha$ and $\gamma^{-1}$ and inverting $\beta$.

\section{Proof of Theorem~1.1}

Let $\overline F_p$ denote the algebraic closure of the field $F_p$ of order $p$, where $p$ is prime.  This is the union of the finite fields $F_q$ for all powers $q$ of $p$, with the natural inclusions, so the group $\overline L:=L_2(\overline F_p)$ is the union of the corresponding groups $L_2(F_q)=L_2(q)$, with the induced inclusions. It follows that any epimorphism $\Delta\to L_2(q)$ can be regarded as a homomorphism $\Delta\to\overline L$, by composition with the natural embedding $L_2(q)\to\overline L$. Conversely, since $\Delta$ is finitely generated, the image of any homomorphism $\Delta\to\overline L$ is contained in a subgroup $L_2(F)$ for some finite subfield $F$ of $\overline F_p$; since $\Delta$ is perfect, so is its image (if nontrivial), so it must be isomorphic to $A_5$ or to $L_2(q)$ for some power $q$ of $p$, by Dickson's classification of the subgroups of $L_2(F)$~\cite[Ch.~XII]{Dic}. We can therefore find the normal subgroups of $\Delta$ with quotient $L_2(q)$ as the kernels of the nontrivial homomorphisms $\theta:\Delta\to\overline L$, excluding any $\theta$ with $\theta(\Delta)\cong A_5$ when $q\neq 4,5$ since $L_2(q)\cong A_5$ if and only $q=4$ or $5$. We find homomorphisms $\theta:\Delta\to\overline L$ by finding elements $\overline\alpha=\theta(\alpha)$, $\overline\beta=\theta(\beta)$ and $\overline\gamma=\theta(\gamma)$ of $\overline L$ which satisfy the defining relations of $\Delta$ in (2.1), using the fact that a non-identity element of $\overline L$ has order $2$, $3$ or $5$ if and only if it has trace $0$, $\pm 1$ or $\pm(-1\pm\sqrt 5)/2$.

\subsection{The restriction to $\Delta_0$}

Assume first that $p\neq 2, 5$. If $\theta:\Delta\to\overline L$ is any nontrivial homomorphism, then its restriction $\psi$ to the simple group $\Delta_0$ must be an isomorphism with a subgroup $G\cong A_5$ of $\overline L$. There is a single conjugacy class of such subgroups $G$ in $\overline L$. This is because, being finite, any pair of such subgroups are both contained in a subgroup $L_2(F)$ for some finite subfield $F$ of $\overline F_p$; now $L_2(F)$ has two conjugacy classes of such subgroups, and these are all conjugate in the subgroup $PGL_2(F)\leq L_2(\widetilde F)\leq\overline L$, where $\widetilde F$ is the quadratic extension of $F$ in $\overline F_p$. Let us define two embeddings $\Delta_0\to\overline L$ to be equivalent if they differ by an inner automorphism of $\overline L$. Given any pair of embeddings $\psi_i\;(i=1,2)$, the conjugacy of their images implies that $\psi_2$ is equivalent to an embedding $\psi_2'$ with the same image as $\psi_1$, so $\psi_2'$ differs from $\psi_1$ by an automorphism of $A_5$. Since $|{\rm Out}\,A_5|=2$ it follows that there are at most two equivalence classes of embeddings of $\Delta_0$. The outer automorphism of $A_5$ transposes its two conjugacy classes of elements of order $5$, and these are distinguished by their images in $\overline L$ having distinct traces $\pm(-1+\sqrt 5)/2$ or $\pm(-1-\sqrt 5)/2$; since traces are invariant under conjugation, it follows that there are exactly two equivalence classes of embeddings $\Delta_0\to\overline L$.

We can construct representatives of these two classes as follows. Let us define $F=F_p$ or $F_{p^2}$ as $p\equiv\pm 1$ or $\pm 2$ mod~$(5)$, so $F$ is the smallest subfield of $\overline F_p$ containing a square root of $5$, or equivalently for which the group $L:=L_2(F)$ contains elements of order $5$. A simple counting argument shows that in any finite field, each element is a sum of two squares, so for each $t=t_i=(-1\pm\sqrt 5)/2\;(i=1,2)$ in $F$ we can find $e,f\in F$ (depending on $i$) which satisfy $e^2+f^2+3=t^2$, that is,
\[e^2+f^2+t+2=0\]
since $t^2=1-t$. We cannot have $e=f=0$, since this
gives $t=-2$ and so $4=t^2=1-t=3$, which is impossible; without loss of generality we may therefore assume that $f\neq 0$. We now define
$\psi_i:\Delta_0\to L\leq\overline L$, for $i=1,2$, by taking $t=t_i$ and sending
\[\alpha\mapsto\overline\alpha:=\frac{1}{2}
\left(\begin{array}{cc}
1-e & -t-f \\
t-f & 1+e
\end{array}\right)
\quad{\rm and}\quad
\beta\mapsto\overline\beta:=\frac{1}{2}
\left(\begin{array}{cc}
t+f & -1-e \\
1-e & t-f
\end{array}\right),
\eqno(3.1)
\]
so that
\[\alpha\beta\mapsto
\left(\begin{array}{cc}
0 & -1 \\
1 & 0
\end{array}\right).\]
Since $\overline\alpha, \overline\beta$ and $\overline\alpha\overline\beta$ have traces $\pm 1, \pm t$ and $0$ they have orders $3$, $5$ and $2$, so it follows that each $\psi_i$ is a homomorphism and hence an embedding by the simplicity of $\Delta_0$. Since the image $\overline\beta$ of $\beta$ under $\psi_i$ has trace $\pm t_i$, the embeddings $\psi_1$ and $\psi_2$ are not equivalent, so every embedding of $\Delta_0$ is equivalent to precisely one of them.

In the case $p=5$ there is a single equivalence class of embeddings $\Delta_0\to\overline L$: the elements of order $5$ in $\Delta_0$ are all represented by elements with trace $\pm 2$, and the outer automorphism of $\Delta_0$ corresponds to conjugacy within the subgroup $PGL_2(5)\cong S_5$ of $\overline L$. We can therefore define $\overline\alpha$ and $\overline\beta$ as above, with $F=F_5$, $e=0$, $f=1$ and $t=2$.

If $p=2$ then $L_2(2^n)$ has a subgroup $G\cong A_5$ if and only if $n$ is even, in which case all such subgroups are conjugate to $L_2(2^2)$. Thus $\overline L$ has a single conjugacy class of such subgroups. There are two equivalence classes of embeddings $\Delta_0\rightarrow\overline{L}$, distinguished by the traces of $\overline\beta$. These are the elements $t_i\in F_4\setminus F_2$, that is, the primitive cube roots of $1$. We define representative embeddings
$\psi_i:\Delta_0\to\overline L$, for $i=1,2$, by taking $t=t_i$ and sending
\[\alpha\mapsto
\left(\begin{array}{cc}
t & 0 \\
t & t +1
\end{array}\right)
\quad{\rm and}\quad
\beta\mapsto
\left(\begin{array}{cc}
0 & t+1 \\
t & t
\end{array}\right),
\eqno(3.2)\]
so
\[\alpha\beta\mapsto
\left(\begin{array}{cc}
0 &1 \\
1 & 0
\end{array}\right).\]

\subsection{Extension to $\Delta$}

A mapping
\[\gamma\mapsto\overline\gamma:=
\left(\begin{array}{cc}
w & x \\
y & z
\end{array}\right)
\in\overline L,\eqno(3.3)\]
with $wz-xy=1$, extends $\psi=\psi_i$ to a homomorphism $\theta=\theta_i:\Delta\to\overline L$ if and only if it preserves the relations
\[\gamma^3=(\beta\gamma)^2=(\alpha\beta\gamma)^2=1\eqno(3.4)\]
of $\Delta$.

Suppose first that $p\neq 2$. Multiplying the matrix in $(3.3)$ by $-1$ if necessary, using $(3.1)$ and taking traces, we can rewrite $(3.4)$ as
\[w+z=1\quad{\rm and }\quad
(t+f)w-(1+e)y+(1-e)x+(t-f)z=y-x=0,\]
or equivalently
\[z=1-w,\quad 2fw=2ex+f-t\quad{\rm and}\quad y=x.\eqno(3.5)\]
The equation $\det\overline\gamma=wz-xy=1$ then becomes
\[4f^2=2f(2ex+f-t)-(2ex+f-t)^2-4f^2x^2,\]
which is equivalent to
\[(e^2+f^2)x^2-etx+\frac{1}{4}(3f^2+t^2)=0.\eqno{(3.6)}\]
This is a quadratic equation for $x$, with coefficients in $F$; by $(3.5)$, extensions of $\psi_i$ to $\Delta$ correspond bijectively to its roots. The discriminant of $(3.6)$ is
\begin{eqnarray*}
D&=&e^2t^2-(e^2+f^2)(3f^2+t^2)\\
&=&-f^2(3e^2+3f^2+t^2)\\
&=&f^2(4t+5)\\
&=&f^2(3\pm 2\sqrt 5)
\end{eqnarray*}
as $t=(-1\pm\sqrt 5)/2$. Since $f\neq 0$, $D$ is a square in $F$ if and only if $3\pm 2\sqrt 5$ is a square, and it is zero if and only if $p=11$ and $t=7$.

We will deal first with the case $p=11$, so that $F=F_{11}$, $\sqrt 5=\pm 4$ and $t=3$ or $7$. If $t=7$ then $D=0$ and we obtain a unique solution $x\in F_{11}$ of $(3.6)$, so $(3.5)$ gives a unique extension of $\psi$ to an epimorphism $\theta:\Delta\to L_2(11)$. Its kernel $K=K_0$, and the associated quotient manifold and tessellation, are studied in detail in~\cite{JM}. If we take $t=4$ then $4t+5=-1$ is a non-square in $F_{11}$, so we obtain two solutions $x$ of $(3.6)$ in the quadratic extension $F(\sqrt D)=F_{11^2}$, giving two extensions of $\psi_i$ to epimorphisms $\Delta\to L_2(11^2)$. The Galois group ${\rm Gal}\,F_{11^2}/F_{11}\cong C_2$ transposes these two roots, and it induces an automorphism of $L_2(11^2)$ fixing $\overline\alpha$ and $\overline\beta$ and transposing the two possible images of $\gamma$, so these two epimorphisms have the same kernel. We therefore have a single normal subgroup of $\Delta$ with quotient $L_2(11^2)$. There are no normal subgroups $K$ with $\Delta/K\cong L_2(11^n)$ for any $n>2$, so this proves Theorem~1.1(c).

From now on we assume that $p\neq 11$, so $D\neq 0$. First suppose that $p\equiv \pm 1$ mod~$(5)$, so that $F=F_p$. Now
\[(3+2\sqrt 5)(3-2\sqrt 5)=-11,\]
and quadratic reciprocity implies that $-11$ is a square mod~$(p)$ if and only if $p$ is a square mod~$(11)$, that is, $p\equiv 1, 3, 4, 5$ or $9$ mod~$(11)$. It follows that for these primes, the two values of $D$, arising from the two choices $t_i$ for $t$, are either both squares in $F_p$ or both non-squares, whereas if $p\equiv 2, 6, 7, 8$ or $10$ mod~$(11)$ then exactly one of them is a square. For each $D$ which is a square we obtain two roots $x\in F_p$ of $(3.6)$, and hence two epimorphisms $\Delta\to L_2(p)$; these have distinct kernels since only the identity element of ${\rm Aut}\,L_2(p)=PGL_2(p)$ fixes the image of $\Delta_0$. Moreover, if both values of $D$ are squares then the resulting four kernels are distinct since $PGL_2(p)$ preserves the trace of $\overline\beta$. If $D$ is a non-square we obtain two roots $x\in F(\sqrt D)=F_{p^2}$, and as in the case $p=11$ these correspond to a single kernel; again, if both values of $D$ are non-squares, the corresponding two kernels are distinct since ${\rm Aut}\,L_2(p^2)=P\Gamma L_2(p^2)$ preserves the trace of $\overline\beta$. To summarise, if $p\equiv 1, 3, 4, 5$ or $9$ mod~$(p)$, as in Theorem~1.1(d), we obtain either four normal subgroups with quotient $L_2(p)$ or two with quotient $L_2(p^2)$, as $3\pm 2\sqrt 5$ are both squares or both non-squares in $F_p$, whereas if $p\equiv 2, 6, 7, 8$ or $9$ mod~$(11)$, as in Theorem~1.1(e), we obtain two normal subgroups with quotient $L_2(p)$ and one with quotient $L_2(p^2)$.

A similar argument applies when $p=5$, except that there is now a single equivalence class of embeddings $\Delta_0\to \overline L$. Since $4t+5=3$ is a non-square in $F=F_5$, ${\rm Gal}\,F_{5^2}/F_5$ transposes the two possible extensions of $\psi$ to $\Delta$, so we obtain a single normal subgroup with quotient $L_2(5^2)$, as in Theorem~1.1(b).

Now suppose that $p\equiv \pm 2$ mod~$(5)$, so that $F=F_{p^2}$. In this case, $-11$ is always a square in $F$, since $F_p(\sqrt{-11})$ is the unique quadratic extension $F_{p^2}$ of $F_p$, so the two discriminants $D$ are either both squares or both non-squares in $F$, as $3\pm 2\sqrt 5$ are or are not both squares. If they are both squares, each yields two epimorphisms $\Delta\to L_2(p^2)$; the resulting four epimorphisms are equivalent in pairs under ${\rm Gal}\,F_{p^2}/F_p$, so each pair have the same kernel, but two epimorphisms corresponding to the same value of $D$ cannot have the same kernel since the centraliser in $P\Gamma L_2(p^2)$ of an icosahedral subgroup is trivial, so there are two normal subgroups with quotient $L_2(p^2)$, as in Theorem~1.1(f). If both values of $D$ are non-squares in $F$, then by considering a further quadratic extension field $F(\sqrt D)=F_{p^4}$ we see that each value of $D$ gives rise to two epimorphisms $\Delta\to L_2(p^4)$; then ${\rm Gal}\,F_{p^4}/F_p\cong C_4$ transposes the two values of $D$, and its subgroup ${\rm Gal}\,F_{p^4}/F_{p^2}\cong C_2$ fixing each $D$ tranposes the two epimorphisms corresponding to $D$, so all four epimorphisms are equivalent under ${\rm Aut}\, L_2(p^4)=P\Gamma L_2(p^4)$ and we obtain a single normal subgroup with quotient $L_2(p^4)$, as in Theorem~1.1(g).

Finally, let $p=2$, so $\psi_i:\Delta_0\to\overline L$ is given by $(3.2)$ and $(3.3)$. The relations $(3.4)$ now give
\[w+z=1,\quad (t+1)y+tx+tz=0\quad{\rm and}\quad x+y=0,\]
leading to
\[w=(t+1)x+1\quad y=x\quad{\rm and}\quad z=(t+1)x.\eqno(3.7)\]
Thus
\[1=\det\overline\gamma=\bigl((t+1)x+1\bigr)(t+1)x-x^2=(t+1)(x^2+x),\]
giving
\[x^2+x+t=0,\eqno(3.8)\]
an irreducible quadratic polynomial for $x$ with coefficients in $F_4$. For each of the two choices for $t$ this polynomial  has two roots $x\in F_{16}$, so we obtain four epimorphisms $\Delta\to L_2(2^4)$. The Galois group ${\rm Gal}\,F_4/F_2$ transposes the two possible values of $t$, and for each $t$ the Galois group ${\rm Gal}\,F_{16}/F_4$ transposes the two roots $x$, so we obtain a single normal subgroup with quotient $L_2(2^4)$, as in Theorem~1.1(a). \hfill$\square$

\section{Normality in $\Omega^+$: proof of Theorem~1.2}

The group $\Omega^+$ is a semidirect product of $\Delta$ by $\langle r\rangle \cong C_2$, with the action of $r$ by conjugation on $\Delta$ given by transposing $\alpha$ and $\gamma^{-1}$ and inverting $\beta$. To prove that each $K$ in Theorem~1.1 is normal in $\Omega^+$ we therefore need to show that there is an automorphism $g=\overline r\in{\rm Aut}\,L_2(q)=P\Gamma L_2(q)$ of $\Delta/K\cong L_2(q)$ acting in the same way on the images $\overline\alpha$, $\overline\beta$ and $\overline\gamma$ of $\alpha$, $\beta$ and $\gamma$; in particular, if $g\in L_2(q)$ we have $\Omega^+/K\cong L_2(q)\times C_2$,  and if $g\in PGL_2(q)\setminus L_2(q)$ we have $\Omega^+/K\cong PGL_2(q)$.

Let $\Pi$ be the set of ordered pairs of elements of order $3$ and $5$ in $L_2(q)$ such that their product has order $2$. Any pair in $\Pi$ generates an icosahedral subgroup $S\cong A_5$ of $L_2(q)$, and any such subgroup contains $120$ pairs, with $S$ permuting them by conjugation  in two regular orbits, and ${\rm Aut}\,S\cong S_5$ permuting them regularly. These icosahedral subgroups are all conjugate in $PGL_2(q)$, with each $S$ equal to its normaliser if $p\neq 5$. In this case it follows that $PGL_2(q)$, acting by conjugation on $\Pi$, has two regular orbits, distinguished by the trace of the element of order $5$ in each pair. If $p=5$ there are one or two conjugacy classes of icosahedral subgroups as $q$ is an odd or even power of $5$; they form a single class in $PGL_2(q)$, each with normaliser $PGL_2(5)\cong S_5$, so in this case $PGL_2(q)$ acts regularly on $\Pi$ (and the elements of order $5$ all have trace $\pm 2$).

Since $(\overline\alpha,\overline\beta)$ and $(\overline\gamma^{-1},\overline\beta^{-1})$ are both pairs in $\Pi$, with $\overline\beta$ and $\overline\beta^{-1}$ having the same trace, the preceding argument shows that there is a unique element $g\in PGL_2(q)$ satisfying $\overline\alpha^g=\overline\gamma^{-1}$ and $\overline\beta^g=\overline\beta^{-1}$. If $p\neq 5$ then the normaliser of $\langle\overline\beta\rangle$ in $PGL_2(q)$ is a dihedral group of order $2(q\pm 1)$ as $5$ divides $q\pm 1$, and the elements inverting $\overline\beta$ are all involutions, so $g^2=1$ and hence $(\overline\gamma^{-1})^g=\overline\alpha$ as required. If $p=5$ then a similar argument applies, except that the normaliser of $\langle\overline\beta\rangle$ is now a Frobenius group, an extension of an elementary abelian Sylow $5$-subgroup of $L_2(q)$ by $C_4$. \hfill$\square$

\medskip

For any given $K$, one can decide the structure of $\Omega^+/K$, and hence that of $N(K)/K$, as follows. The element $g\in PGL_2(q)$ in the proof of Theorem~B is in $L_2(q)$ or not as $\det(g)$ is or is not a square in $F_q$, giving $\Omega^+/K\cong L_2(q)\times C_2$ or $PGL_2(q)$ respectively. The necessary and sufficient conditions
\[g^2=1,\quad \overline\beta^g=\overline\beta^{-1},\quad \overline\alpha^g=\overline\gamma^{-1}\]
on $g$ can be expressed as linear equations in the entries $g_i$ of the matrix
\[G=
\left(\begin{array}{cc}
g_1 & g_2 \\
g_3 & g_4
\end{array}\right)
\in GL_2(q)\]
chosen to represent $g$. Thus the condition $g^2=1$ gives
\[g_1+g_4=0,\]
in which case $\overline\beta^g=\overline\beta^{-1}$ can be written as $(\overline\beta g)^2=1$. If $p\neq 2$ then by $(3.1)$ this gives
\[2fg_1+(1-e)g_2-(1+e)g_3=0.\]
If $A$ and $C$ are the matrices representing $\overline\alpha$ and $\overline\gamma$ in $(3.1)$ and $(3.3)$, then the condition $\overline\alpha^g=\overline\gamma^{-1}$ can be written as $G^{-1}AG=\lambda C^{-1}$ where $\lambda\in F_q\setminus\{0\}$. Since $A$ and $C$, and hence $G^{-1}AG$ and $C$, have trace $1$, it follows that $\lambda=1$, so $G^{-1}AG=C^{-1}$. Writing this as $AG=GC^{-1}$, putting $g_4=-g_1$, and comparing the entries on each side, we obtain four more equations
\begin{eqnarray*}
(2ex-fe-t)g_1+2fxg_2-(t+f)fg_3&=&0,\\
f(f+2x+t)g_1+(-fe+t-2xe)g_2&=&0,\\
f(t-f-2x)g_1+(-t+fe+2xe)g_3&=&0,\\
(-fe+2ex-t)g_1+(t-f)fg_2+2fxg_3&=&0.
\end{eqnarray*}

The uniqueness of $g$ in $PGL_2(q)$ shows that these six linear equations determine a unique point $[g_1, g_2, g_3, g_4]$ in the subset $\det(g)=g_1g_4-g_2g_3\neq 0$ of projective $3$-space $P^3(q)$ over $F_q$. It follows that they determine a unique point in $P^3(q)$, for otherwise there would be a projective line of such points, containing at least two points not on the quadratic curve $\det(g)=0$ since $q\geq 3$. Thus the solution vector $(g_1, g_2, g_3, g_4)$ is unique, up to scalar multiplication.

Suppose first that $f+t+x\neq 0$. Since $f\neq 0$ we have
\[g_1=\frac{-t+ef+2ex}{ f(f+t+2x)}g_2\]
and
\[g_3=\frac{f-t+2x}{t+f+2x},\]
so that
\[\det(g)=-g_1^2-g_2g_3=\frac{g_2^2 s}{(f+t+2x)^2}\]
where
\[s=f^{-2}\big((8+4t)x^2+(4et+8f+4ft)x-1+t+2eft+3f^2\bigr).\]
Thus $\det(g)$ is a square if and only if $s$ is a square. Now for a given $t$ the two possible roots $x$ of the quadratic equation $(3.6)$ are
\[x={et\pm\sqrt D\over2(e^2+f^2)}={et\pm f\sqrt{4t+5}\over-2(t+2)},\]
so a calculation with MAPLE yields
\[s=-2(-3\pm\sqrt{4t+5})\]
as we choose the $+$ or $-$ sign in the formula for $x$. Thus $g\in L_2(q)$ or not, giving $\Omega^+/K\cong L_2(q)\times C_2$ or $PGL_2(q)$ respectively, as the element
\[s=-2(-3\pm\sqrt{4t+5})=6\mp2\sqrt{3\pm 2\sqrt 5}\]
is or is not a square in $F_q$.

If $f+t+x=0$ then $t-f-x\neq 0$ and we can use a similar argument, starting with
\[g_1=\frac{2ex+ef-t}{f(2x+f-t)}g_3\]
and leading to the same criterion.

Note that for a given value of $t=(-1\pm\sqrt 5)/2$, the product of the two possible values of $s$ is
\[(-2(-3+\sqrt{4t+5}))(-2(-3-\sqrt{4t+5}))=(4t)^2,\]
so either both are squares in $F_q$, or both are non-squares; it follows that even if the two roots $x$ determine distinct normal subgroups $K$, the quotients $\Omega^+/K$ corresponding to a given value of $t$ are isomorphic. However, as shown by Example~4.3 below, quotients $\Omega^+/K$ for different values of $t$ need not be isomorphic, even when the quotients $\Delta/K$ are isomorphic.

The case $p=2$ can be dealt with in a similar way, using $(3.2)$ instead of $(3.1)$, but we will omit the details since the result follows immediately from Torstensson's work. She has shown in~\cite{Tor} that $\Omega^+$ (denoted there by $\Gamma$) has a normal subgroup $N$ with $\Omega^+/N\cong L_2(2^4)$; then $N\cap\Delta$ is the subgroup $K$ in Theorem~1.1(a), and this is normal in $\Omega^+$ with quotient $L_2(2^4)\times C_2$.

The criterion given above confirms Lemmas~9, 11 and 14 of~\cite{JM}, which show by individual calculations that $\Omega^+/K\cong PGL_2(q)$ when $q=3^2$, $11$ or $19$.

\medskip

\noindent{\bf Example 4.1.} If $p=3$ then $3\pm 2\sqrt 5=4t+5=t-1=(t+1)^2$ are both squares in $F=F(t)=F_9$, so Theorem~1.1(f) implies that there are two normal subgroups $K$ of $\Delta$ with $\Delta/K\cong L_2(3^2)$ (these are the groups $K_1$ and $K_2$ studied in~\cite{JM}). In this case $s=6\mp 2(t+1)=\pm(t+1)$, both of which are non-squares in $F_9$, so $\Omega^+/K\cong PGL_2(3^2)$ for each $K$.

\medskip

\noindent{\bf Example 4.2.} If $p=5$ then Theorem~1.1(b) gives a unique normal subgroup $K$ of $\Delta$ with $\Delta/K\cong L_2(5^2)$. In this case $s=1\mp 2\sqrt 3=(2\pm 2\sqrt 3)^2$ in $F_{25}$, so $\Omega^+/K\cong L_2(5^2)\times C_2$.

\medskip

\noindent{\bf Example 4.3.} If $p=59$ then $\sqrt 5=\pm 8$ in $F_p$, so $3\pm 2\sqrt 5$ takes the value $19=14^2$ or $-13=20^2$. Theorem~1.1(d) therefore shows that there are four normal subgroups $K$ of $\Delta$ with $\Delta/K\cong L_2(59)$. Taking $3\pm 2\sqrt 5=19$ gives $s=6\mp 28=-22$ or $-25$, both non-squares in $F_{59}$ by quadratic reciprocity, so $\Omega^+/K\cong PGL_2(59)$ for the corresponding pair of subgroups $K$. However, taking $3\pm 2\sqrt 5=-13$ gives $s=25=5^2$ or $-13=20^2$, so $\Omega^+/K\cong L_2(59)\times C_2$ for the other pair.

\medskip

In \S 11 we give tables of results for primes $p\leq 251$, and in each case the structure of $\Omega^+/K$ is found by using GAP to check whether or not $s$ is a square in $F_q$. It is striking that in all cases in these tables where $F_q>F$ (or equivalently $3\pm 2\sqrt 5$ is not a square  in $F$), we have $\Omega^+/K\cong L_2(q)\times C_2$. On the basis of this evidence we conjecture that $\Omega^+/K\cong L_2(q)\times C_2$ whenever $F_q>F$. A similar phenomenon for the tetrahedral group $[5,3,5]^+$ is observed in~\cite{JL}, where the analogues of Theorems~1.1 and 1.2 are proved for this group.

\section{Normality in $\Gamma$: proof of Theorem 1.3}

Let $K$ be a normal subgroup of $\Delta$ with $\Delta/K\cong L_2(q)$; we need to determine whether $K$ is normal in $\Gamma$, and if so to determine the structure of $\Gamma/K$.

The element $i=(abc)^5$ is the central involution in the Coxeter subgroup $\Gamma_0\cong A_5\times C_2$ of $\Gamma$. Since $\Gamma=\langle\Delta,i\rangle$ it follows that $K$ is normal in $\Gamma$ if and only if $K^i=K$, or equivalently $L_2(q)$ has an automorphism $\iota$ corresponding to that induced on $\Delta$ by conjugation by $i$; in this case $\Gamma/K$ is a semidirect product of $L_2(q)$ by $C_2$, with the latter acting by conjugation on $L_2(q)$ as $\langle\iota\rangle$. Now $i$ commutes with $\alpha$ and $\beta$, so if $\iota$ exists then it must fix the images $\overline\alpha$ and $\overline\beta$ of $\alpha$ and $\beta$ in $L_2(q)$. Conjugating the relations
\[\gamma^3=(\beta\gamma)^2=(\alpha\beta\gamma)^2=1\]
of $\Delta$ by $i$, we see that the element $\delta:=\gamma^i$ satisfies similar relations
\[\delta^3=(\beta\delta)^2=(\alpha\beta\delta)^2=1.\]
This means that the entries of a matrix representing the image $\overline\delta$ of $\delta$ must satisfy the same polynomial equations over $F$ (given in \S 3.2) as the corresponding entries $w, x, y$ and $z$ for $\overline\gamma$. We have seen that these equations have at most two solutions, determined by the choice of a root $x$ of a quadratic polynomial $Q(x)$ (see equations~$(3.6)$ and $(3.8)$ for $p\neq 2$ and $p=2$), giving at most two possibilities for $\overline\gamma$ and $\overline\delta$. There are now three cases to consider:

\begin{enumerate}

\item If $\overline\gamma=\overline\delta$ then $\iota$ is the identity automorphism, so $K$ is normal in $\Gamma$ with $\Gamma/K\cong L_2(q)\times C_2$.

\item If $\overline\gamma\neq\overline\delta$ and $Q(x)$ is irreducible over $F$, then $F_q$ is a quadratic extension $F(x)$ of $F$; the Galois group of this extension transposes the two solutions and induces the required automorphism $\iota$ of $L_2(q)$, fixing $\overline\alpha$ and $\overline\beta$ and transposing $\overline\gamma$ and $\overline\delta$; thus $K$ is normal in $\Gamma$ with quotient $P\Sigma L_2(q)$ or $P\Sigma L_2(q)^+$ as $q=p^2$ or $p^4$.

\item If $\overline\gamma\neq\overline\delta$ and $Q(x)$ is reducible over $F$, then $F_q=F$; in this case there is no such automorphism $\iota$ since only the identity automorphism fixes an icosahedral subgroup of $L_2(q)$, so $K$ is not normal in $\Gamma$.

\end{enumerate}

We saw in \S 3.2 that if $p=11$ and $t=7$ then the polynomial $Q(x)$ has discriminant $D=0$, so it has one repeated root $x$; thus $\overline\gamma=\overline\delta$, so the corresponding group $K=K_0$, the first group in Theorem~1.1(c), is normal in $\Gamma$ with quotient $L_2(11)\times C_2$. In fact, Theorem~5 of~\cite{Pao} shows that $L_2(q)\times C_2$ is a quotient of $\Gamma$ (denoted there by $C_{5,2}$) only when $q=11$, and Theorem~1.1 shows that $K_0$ is the unique normal subgroup with such a quotient. For all other $K$ either (2) or (3) must occur.

If $p\neq 2$ then $Q(x)$, given by equation~$(3.6)$, is reducible if and only if its discriminant $D$ is a square in $F$. Thus (3) occurs if and only if the element $4t+5=3\pm 2\sqrt 5$ is a non-zero square in $F$, as happens for instance when $q=11^2$ and $t=3$. This deals with all cases except case~(a), where $p=2$.

If $p=2$ then Theorem~1.1(a) shows that there is one subgroup $K$ with $q=2^4$. As shown in equation $(3.8)$, this leads to the polynomial $Q(x)=x^2+x+t$ which is irreducible over $F=F_4$, so $K$ is normal in $\Gamma$ with quotient $P\Sigma L_2(2^4)^+$. \hfill$\square$

\section{Determining $N(K)$: proof of Corollary 1.4}

Since $N(K)$ is a discrete subgroup of ${\rm Iso}\,{\cal H}^3$ containing $\Delta$, results of Derevnin and Mednykh~\cite{DM} show that $N(K)\leq\Omega$. By Theorem~1.2 we have $\Omega^+\leq N(K)$. Since $|\Omega:\Omega^+|=2$ it follows that $N(K)=\Omega$ or $\Omega^+$, and since $\Omega$ is generated by $\Gamma$ and $\Omega^+$ these two cases correspond to $K$ being normal or not normal in $\Gamma$.\hfill$\square$

\section{Proof of Corollary 1.5}

Let $N$ be a proper normal subgroup of $\Gamma$ with $\Gamma/N\cong S\leq L_2(q)$ for some prime power $q$. If $N$ is a proper subgroup of $\Delta$ then $N$ is torsion-free (see \S 2), so $\Gamma_0$ is mapped isomorphically into $\Gamma/N$ and hence into $L_2(q)$. This is impossible, since $\Gamma_0$ is nonsolvable, with a central involution, whereas the centraliser of any involution in $L_2(q)$ is solvable.

We may therefore assume that $N\not\leq\Delta$, so that $N^+:=N\cap\Delta$ is a normal subgroup of $\Gamma$ with $\Delta/N^+\cong S$ and $\Gamma/N^+\cong S\times C_2$. Since $\Delta$ is perfect, so is $S$, so $S\cong L_2(q')$ for some prime power $q'$ by~\cite[Ch.~XII]{Dic}. The normal subgroups $K$ of $\Delta$ with such a quotient are classified in Theorem~1.1, and by Theorem~1.3 the only instance in which $K$ is normal in $\Gamma$ with quotient $L_2(q')\times C_2$ is case~(c), with $q'=11$ and $K=K_0$. It follows that $N^+=K_0$, so $N=N_0$, the only normal subgroup of $\Gamma$ containing $K_0$ with index $2$, and $\Gamma/N\cong\Delta/K_0\cong L_2(11)$.\hfill$\square$

\section{Action of $L_2(q)$ on quotient tessellations ${\cal I}/K$}

If $K$ is a normal subgroup of $\Delta$ with quotient $L_2(q)$ then the manifold ${\cal M}={\cal H}^3/K$ has an icosahedral tessellation ${\cal I}/K$, with $L_2(q)$ acting as its orientation-preserving symmetry group $\Delta/K$. If $K$ is normal in $\Gamma$ then the full symmetry group is $\Gamma/K$, described in Theorem~1.3, including orientation-reversing symmetries, whereas if $K$ is not normal in $\Gamma$ then ${\cal I}/K$ is one of a chiral pair of tessellations, with only orientation-preserving symmetries. Since the elements of $\Omega\setminus\Gamma$ transform $\cal I$ to its dual tessellation ${\cal I}^*\cong{\cal I}$, the fact that $K$ is normal in $\Omega^+$ (see Theorem~1.2) implies that each element of $\Omega^+\setminus\Delta$ induces an isometry of $\cal M$ which realises an isomorphism between ${\cal I}/K$ and its dual. In particular this applies to the half-turn $r\in\Omega^+\setminus\Delta$.

\subsection{Actions on cells and vertices}

In the action of $L_2(q)$ on ${\cal I}/K$, the cells and the vertices are each permuted transitively, and their stabilisers are the images in $L_2(q)$ of $\Delta_0=\langle\alpha, \beta\rangle$ and of $\Delta_0^r=\langle\beta, \gamma\rangle$. These icosahedral subgroups $I$ of $L_2(q)$ have index
\[|L_2(q):I|=\frac{|L_2(q)|}{60}=\left\{
\begin{array}{ll}
q(q^2-1)/120 & \mbox{if $q$ is odd,}\\
& \\
68 & \mbox{if $q=2^4$,}
\end{array}
\right.
\eqno(8.1)\]
so this is the number of cells and also the number of vertices.

In order to understand further the actions of $L_2(q)$ on cells and on vertices, it  is useful to summarise here the existence and conjugacy of its icosahedral subgroups~\cite[Ch.~XII]{Dic}. These exist if and only if $q\equiv 0$ or $\pm 1$ mod~$(5)$. If $p\neq 2$ or $5$, or if $q=5^n$ with $n$ even, then they form two conjugacy classes, which merge to form a single conjugacy class in $PGL_2(q)$. There is a single class of icosahedral subgroups in $L_2(q)$ if $q=2^n$ with $n$ even, or if $q=5^n$ with $n$ odd, and there are none if $q=2^n$ with $n$ odd. These subgroups are all self-normalising in $L_2(q)$, except when $q=5^n$ with $n$ even, in which case they have normaliser isomorphic to $S_5\cong PGL_2(5)$.

It follows that if $L_2(q)$ is a quotient of $\Delta$, as in Theorem~1.1, then except in case~(a), where $q=2^4$, it has two conjugacy classes of icosahedral subgroups, which are conjugate to each other in $PGL_2(q)$. The image $\overline r$ of $r$ transforms ${\cal I}/K$ to its dual tessellation, so if $\Omega^+/K\cong PGL_2(q)$ then the stabilisers of cells and of vertices form distinct conjugacy classes in $L_2(q)$, and hence these two actions of $L_2(q)$ are inequivalent; if $\Omega^+/K\cong L_2(q)\times C_2$, however, the icosahedral subgroups in one class are stabilisers of both cells and vertices, while those in the other class stabilise neither, so the two actions are equivalent. In case~(a) there is a single conjugacy class of icosahedral subgroups of $L_2(2^4)$, stabilising both cells and vertices, so again the two actions are equivalent.

Except in case~(b), where $q=5^2$, each stabiliser $I$ is its own normaliser in $L_2(q)$, so it stabilises a single cell or vertex. In case~(b), however, each stabiliser has index $2$ in its normaliser $N(I)\cong S_5\cong PGL_2(5)$, so it stabilises a pair of cells and, since $\Omega^+/K\cong L_2(5^2)\times C_2$ (see Example~4.2), a pair of vertices. In this latter case, these pairs form blocks of imprimitivity for the actions of $L_2(5^2)$ on cells and on vertices. More generally, if $F_q>F$ then each stabiliser $I$ is contained in proper subgroups isomorphic to $L_2(F)$ and to $PGL_2(F)$, so $L_2(q)$ acts imprimitively on cells and on vertices, whereas if $F_q=F$ then each $I$ is a maximal subgroup of $L_2(q)$, so the cells and vertices are permuted primitively.

One can use this information to compute the cycle structures of elements of $L_2(q)$, acting on cells or vertices. The orders of such elements are $p$ and the divisors of $(q\pm 1)/2$ when $q$ is odd, and the divisors of $2, 15$ and $17$ when $q=2^4$. It is convenient first to consider elements of orders $2, 3$ and $5$, since these are the only non-identity elements which can have fixed points in these actions.

If  $q\equiv\pm 1$ mod~$(4)$ there are $q(q\pm 1)/2$ involutions in $L_2(q)$. The tessellation ${\cal I}/K$ has $q(q^2-1)/120$ vertices, each fixed by the $15$ involutions in its icosahedral stabiliser. Since the involutions in $L_2(q)$ are all conjugate, each fixes
\[\frac{q(q^2-1)}{120}\times 15\times \frac{2}{q(q\pm 1)}=\frac{q\mp 1}{4}\]
vertices, and similarly it leaves invariant the same number of cells; the remaining
\[\frac{q(q^2-1)}{120}-\frac{(q\mp 1)}{4}=\frac{(q\mp 1)\bigl(q(q\pm 1)-30\bigr)}{120}\]
vertices or cells are permuted as $2$-cycles. A similar argument applies if $q=2^4$, except that now there are $q^2-1=255$ involutions, each fixing four of the $q(q^2-1)/60=68$ vertices or cells, and having $(68-4)/2=32$ $2$-cycles on the others.

If $q\equiv\pm 1$ mod~$(3)$ there is a single conjugacy class of $q(q\pm 1)$ elements of order $3$ in $L_2(q)$, with $20$ of them fixing each vertex or cell; in either action, if $q$ is odd an element of order $3$ has $(q\mp 1)/6$ fixed points and $(q\mp 1)\bigl(q(q\pm 1)-20\bigr)/360$ $3$-cycles, while if $q=2^4$ the $272$ elements of order $3$ each have five fixed points and $21$ $3$-cycles. If $q=3^2$ the elements of order $3$ form two conjugacy classes of size $20$; those in one class, represented by $\overline\alpha$, have cycle structures $3^2$ and $1^33$ on the vertices and cells, while those in the other class, represented by $\overline\gamma$, have cycle structures $1^33$ and $3^2$ respectively.

If $q\equiv\pm 1$ mod~$(5)$ there are two conjugacy classes of $q(q\pm 1)$ elements of order $5$ in $L_2(q)$, each consisting of the squares of the elements in the other class. There are $24$ elements of order $5$ fixing each vertex or cell, so if $q$ is odd then each element of order $5$ has $(q\mp 1)/10$ fixed points and $(q\mp 1)\bigl(q(q\pm 1)-12\bigr)/600$ $5$-cycles on the vertices and on the cells, while if $q=2^4$ the $544$ elements of order $5$ each have three fixed points and $13$ $5$-cycles. If $q=5^2$ the elements of order $5$ form two conjugacy classes of size $312$, each closed under squaring (and hence also inversion); these are distinguished by which conjugacy class of icosahedral subgroups they lie in, so the elements in one class, represented by $\overline\beta$, have ten fixed points and $24$ $5$-cycles on the vertices and on the cells, while those in the other class have no fixed points and $26$ $5$-cycles.

One can use these results on elements of orders $2$, $3$ and $5$ to determine the cycle structures of elements of other orders. We will restrict ourselves to a simple example:

\medskip

\noindent{\bf Example 8.1.} Suppose that $q\equiv -1$ mod~$(12)$, and that an element $g\in L_2(q)$ has order~$6$. In its action on vertices, $g$ has no fixed points, since icosahedral subgroups have no elements of order $6$, so $g$ consists of disjoint cycles of lengths $2$, $3$ and $6$. Now $g^2$ has order $3$, so it fixes $(q+1)/6$ vertices, and hence $g$ contains $(q+1)/12$ $2$-cycles. Similarly the involution $g^3$ fixes $(q+1)/4$ vertices, so $g$ contains $(q+1)/12$ $3$-cycles, and the remaining
\[\frac{q(q^2-1)}{120} - \frac{(q+1)}{6} - \frac{(q+1)}{4} = \frac{(q+1)(q^2-q-50)}{120}\]
vertices form $(q+1)(q^2-q-50)/20$ $6$-cycles. The cycle structure of $g$ on cells is the same.

\medskip

Other orders of elements and values of $q$ can be treated in the same way, and this allows one to determine the permutation characters $\pi_0(g)$ and $\pi_3(g)$, giving the numbers of vertices and cells fixed by each $g\in L_2(q)$. In fact, our preceding analysis of fixed points shows that $\pi_0=\pi_3$ in all cases except $q=3^2$ and $5^2$, where these characters differ on elements of orders $3$ and $5$ respectively and are transposed by conjugation in $PGL_2(q)$.

Since these actions of $L_2(q)$ are transitive, each $\pi_i$ is formed by inducing the principal character of a stabiliser $I$ up to $L_2(q)$, so by Frobenius reciprocity [Dor, Theorem~9.4(c)] the multiplicity $n_{\chi}$ of any irreducible complex character $\chi$ of $L_2(q)$ in $\pi_i$ is given by
\[n_{\chi}=\frac{1}{60}\sum_{g\in I}\chi(g)=\frac{\chi(g_1)}{60}+\frac{\chi(g_2)}{4}+\frac{\chi(g_3)}{3}+\frac{\chi(g_5)}{5}+\frac{\chi(g_5^2)}{5},\eqno(8.2)\]
where $g_j$ has order $j$ in $I$. The character tables for the groups $L_2(q)$ are easily obtained from those for $SL_2(q)$ computed by Schur in~\cite{Sch} (see also~\cite[\S 38]{Dor}), and using these one can now find the decompositions of $\pi_0$ and $\pi_3$ as sums of irreducible characters. There are many different cases to consider, depending on the congruence class of $q$ mod~$(60)$, so we will simply give an illustrative example.

\medskip

\noindent{\bf Example 8.2.} Suppose that $q\equiv 11$ mod~$(60)$. Since $q\equiv 3$ mod~$(4)$, $L_2(q)$ has, in addition to the principal character, two complex conjugate irreducible characters of degree $(q-1)/2$, one of degree $q$, and $(q-3)/4$ each of degrees $q-1$ and $q+1$. (These details are slightly different  for $q\equiv 1$ mod~$(4)$ and for $q=2^n$.) Let $g_{\pm}$ be elements of orders $(q\pm 1)/2$, so we can take $g_2=g_+^{(q+1)/4}$, $g_3=g_+^{(q+1)/6}$ and $g_5=g_-^{(q-1)/10}$.

The two irreducible characters $\chi$ of degree $\chi(1)=(q-1)/2$ satisfy $\chi(g_2)=(-1)^{(q-1)/4}$, $\chi(g_3)=(-1)^{(q+7)/6}=-1$ and $\chi(g_5)=\chi(g_5^2)=0$, so by $(8.2)$ they have multiplicity
\[n_{\chi}=\left\{
\begin{array}{ll}
(q-11)/120 & \mbox{if $q\equiv 11$ mod~$(120)$,}\\
(q-71)/120 & \mbox{if $q\equiv 71$ mod~$(120)$.}
\end{array}\right.
\]

The irreducible character $\chi$ of degree $\chi(1)=q$ satisfies $\chi(g_2)=\chi(g_3)=-1$ and $\chi(g_5)=\chi(g_5^2)=1$, so it has multiplicity
$$n_{\chi}=(q-11)/60.$$

The $(q-3)/4$ irreducible characters $\chi$ of degree $\chi(1)=q-1$ satisfy $\chi(g_-^j)=0$ for all $g_-^j\neq 1$, so that $\chi(g_5)=\chi(g_5^2)=0$; they also satisfy $\chi(g_+^j)=-e^{4\pi ijk/(q+1)}-e^{-4\pi ijk/(q+1)}$ for all $g_+^j\neq 1$, where $k=1,\ldots, (q-3)/4$, so $\chi(g_2)=-2(-1)^k$ and $\chi(g_3)=-e^{2\pi ik/3}-e^{-2\pi ik/3}$; they therefore have multiplicity
\[n_{\chi}=\left\{
\begin{array}{ll}
(q-71)/60 & \mbox{if $k\equiv 0$ mod~$(6)$,}\\
(q+49)/60 & \mbox{if $k\equiv\pm 1$ mod~$(6)$,}\\
(q-11)/60 & \mbox{if $k\equiv \pm 2$ or $3$ mod~$(6)$.}
\end{array}\right.\]

The $(q-3)/4$ irreducible characters $\chi$ of degree $\chi(1)=q+1$ satisfy $\chi(g_+^j)=0$ for all $g_+^j\neq 1$, so that $\chi(g_2)=\chi(g_3)=0$; they also satisfy $\chi(g_-^j)=e^{4\pi ijk/(q-1)}+e^{-4\pi ijk/(q-1)}$ for all $g_-^j\neq 1$, where $k=1,\ldots, (q-3)/4$, so that $\chi(g_5)=e^{2\pi ik/5}+e^{-2\pi ik/5}$ and $\chi(g_5^2)=e^{4\pi ik/5}+e^{-4\pi ik/5}$; it follows that they have multiplicity
\[n_{\chi}=\left\{
\begin{array}{ll}
(q+49)/60 & \mbox{if $k\equiv 0$ mod~$(5)$,}\\
(q-11)/60 & \mbox{if $k\not\equiv 0$ mod~$(5)$.}
\end{array}\right.\]

When $q=11$, for instance, this shows that the permutation characters $\pi_0$ and $\pi_3$ of $L_2(11)$ both decompose as $\chi_1+\chi_5$ in ATLAS notation, where $\chi_1$ is the principal character and $\chi_5$ is an irreducible character of degree $10$; the fact that $\pi_i-\chi_1$ is irreducible shows that $L_2(11)$ acts as a doubly transitive group on the vertices and on the cells. Indeed, since $L_2(q)$ has no irreducible characters of degree greater than $q+1$ for any $q$, it easily follows that the only other case where it is doubly transitive on the vertices or cells is when $q=3^2$. In this case $L_2(3^2)$ acts on the six vertices and the six cells as the quadruply transitive group $A_6$, with $\pi_0=\chi_1+\chi_3$ and $\pi_3=\chi_1+\chi_2$ in ATLAS notation if we choose $\overline\alpha$ and $\overline\gamma$ to be in the conjugacy classes 3A and 3B respectively. (Reversing this choice, by an outer automorphism of $A_6$, transposes $\chi_2$ and $\chi_3$.)

\medskip

The cases $q=3^2$ and $q=11$ are the only instances where the rank (number of orbits of a stabiliser) is $2$. More generally, for any $q$ these two actions of $L_2(q)$ have rank
\[\frac{1}{|L_2(q)|}\sum_{g\in L_2(q)}\negthinspace\pi_i(g)^2=\sum_{\chi}n_{\chi}^2.\]

\subsection{Actions on faces and edges}

The faces and the edges of ${\cal I}/K$ are also permuted transitively by $L_2(q)$, and these actions can be analysed in a similar way. The main difference is that the (setwise) stabilisers are now dihedral subgroups of order $6$, isomorphic to $S_3$, so the number of faces and of edges is $|L_2(q)|/6$, which is $10$ times the number of cells and of vertices given in (8.1).

Each face is invariant under a rotation group of order $3$, conjugate to $\langle\overline\alpha\rangle$, which preserves its two incident cells, and also under half-turns about its three axes of symmetry, which transpose these two cells. Similarly each edge is fixed pointwise by a rotation group of order $3$, conjugate to $\langle\overline\gamma\rangle$, which permutes its incident cells in a $3$-cycle, and it is reversed by three half-turns, each of which transposes two of these cells. The actions of $L_2(q)$ on faces and on edges are imprimitive, since these stabilisers are never maximal subgroups of $L_2(q)$: for instance they are contained in dihedral subgroups of order $q\mp 1$ if $q\equiv\pm 1$ mod~$(6)$, or of order $12$ or $30$ if $q=3^2$ or $2^4$. If $q\equiv\pm 5$ mod~$(12)$ or $q=3^2$ or $2^4$ then these stabilisers are self-normalising in $L_2(q)$, so they each stabilise a single face or edge; if $q\equiv\pm 1$ mod~$(12)$ then they have index $2$ in their normalisers (dihedral groups of order $12$), so they each stabilise two faces or edges. Except when $q=3^2$, the subgroups of order $3$ in $L_2(q)$ form a single conjugacy class, so they each stabilise an edge and a face, or a pair of each if $q\equiv\pm 1$ mod~$(12)$. However, there are two conjugacy classes of subgroups of order $3$ in $L_2(3^2)$, represented by $\langle\overline\alpha\rangle$ and $\langle\overline\gamma\rangle$, so these subgroups respectively stabilise one face and no edges, or vice versa.

If $q\equiv\pm 1$ mod~$(4)$ then each involution stabilises
\[\frac{q(q^2-1)}{12}\times 3\times\frac{2}{q(q\pm 1)}=\frac{q\mp 1}{2}\]
faces and edges, but if $q=2^4$ then the number is $680\times 3\div 255=8$. Similarly, if $q\equiv\pm 1$ mod~$(6)$ then each element of order $3$ stabilises $(q\mp 1)/6$ faces and edges, but if $q=2^4$ then the number is $680\times 2\div 272=5$. If $q=3^2$ the $20$ elements of order $3$ conjugate to $\overline\alpha$ stabilise $60\times 2\div 20=6$ faces and no edges, while the $20$ conjugate to $\overline\gamma$ stabilise six edges and no faces.

As in \S 8.1 one can use these results to find the cycle structures of elements of other orders. For instance, an element of order $6$ in $L_2(11)$ has cycle structure $1^02^13^26^{17}$ on both edges and faces.

If $\pi_1$ and $\pi_2$ are the permutation characters of $L_2(q)$ on edges and faces then the above argument shows that $\pi_1=\pi_2$ in all cases except $q=3^2$, when these characters differ on the elements of order $3$.  The multiplicity $n_{\chi}$ of any irreducible complex character $\chi$ of $L_2(q)$ in each $\pi_i$ is equal to
\[\frac{1}{6}\sum_{g\in H}\chi(g)=\frac{\chi(g_1)}{6}+\frac{\chi(g_2)}{2}+\frac{\chi(g_3)}{3},\]
where $g_j$ has order $j$ in a stabiliser $H$. For instance, when $q=11$ we have
\[\pi_1=\pi_2=\chi_1+\chi_2+\chi_3+\chi_4+3\chi_5+\chi_6+2\chi_7+2\chi_8,\]
where $\chi_1,\ldots, \chi_8$ are irreducible characters of $L_2(11)$ of degrees $1, 5, 5, 10, 10, 11, 12$ and $12$ in ATLAS notation, so both actions have rank $\sum_{\chi}n_{\chi}^2=22$.

\section{Action of $L_2(q)$ on quotient manifolds ${\cal M}={\cal H}^3/K$}

If $K$ is a normal subgroup of $\Delta$ with quotient $L_2(q)$ then the manifold ${\cal M}={\cal H}^3/K$ admits an action of $L_2(q)$ as a group of orientation-preserving isometries. This is a subgroup of index $4$ or $2$ in the full isometry group $\Omega/K$ or $\Omega^+/K$ of $\cal M$ as $K$ is or is not normal in $\Gamma$. If a non-identity element $g\in L_2(q)$ has fixed points in $\cal M$ then these form a disjoint family of geodesics, each an axis around which $g$ induces a rotation. Since $\cal M$ is compact these geodesics have finite length, so they are homeomorphic to circles, and there are finitely many of them, so they form a link in $\cal M$. Since $g$ preserves the tessellation ${\cal I}/K$ it leaves at least  one vertex, edge, face or cell invariant, so it has order $2$, $3$ or $5$. In each of these cases we can use the results of \S 8 to study the action of $g$ on $\cal M$.

\subsection{Elements of order $3$}

Suppose that $g$ has order $3$, and that $q\neq 3^2$. The arguments in \S 8.2 then show that $g$ leaves invariant a face, so it has an axis of rotation passing perpendicularly through the centre of that face. This axis must be enclosed by a `bracelet' of $\langle g\rangle$-invariant icosahedral cells, each adjacent pair meeting across a $\langle g\rangle$-invariant face. Since the faces of ${\cal I}/K$ are permuted transitively by $L_2(q)$, each of them lies in such a bracelet for some element of order $3$. The normaliser $N(\langle g\rangle)$ of $\langle g\rangle$ in $L_2(q)$ acts transitively on the $\langle g\rangle$-invariant faces, so the bracelets associated with $g$ all contain the same number of cells; since the elements of order $3$ are all conjugate in $L_2(q)$, this number $n_3$ is the same for all such $g$. We have seen that there are $(q\mp 1)/6$ $\langle g\rangle$-invariant cells if $q\equiv\pm 1$ mod~$(6)$, or five if $q=2^4$, so the number of bracelets for each $g$ is respectively $(q\mp 1)/6n_3$ or $5/n_3$. Since $K$ is normal in $\Omega^+$, a similar argument using duality shows that $g$ also has the same number of axes of rotation consisting of $n_3$ vertices and  $n_3$ edges of ${\cal I}/K$. These two types are the only axes of rotation for $g$; they form two orbits under $N(\langle g\rangle)$, and as $g$ ranges over all the elements of order $3$ the corresponding axes form two orbits under $L_2(q)$, transposed by isometries induced by $\Omega^+\setminus\Delta$. When $q=3^2$, the conjugates of $\overline\alpha$ have $3/n_3$ axes, each enclosed by a bracelet of $n_3$ cells, while the conjugates of $\overline\gamma$ have $3/n_3$ axes, each consisting of $n_3$ vertices and $n_3$ edges; again these two types of axes form two orbits under $L_2(q)$, transposed by $\Omega^+\setminus\Delta$.

In order to evaluate $n_3$ we consider the element $(abc)^5d$ of $\Delta$. This is the composition of the central involution $(abc)^5$ of $\Gamma_0$, inducing the antipodal isometry of an icosahedral cell of $\cal I$, and a reflection $d$ in one of its faces, so it acts on $\cal I$ as a screw transformation, translating this cell to an adjacent cell and rotating by a half-turn so that these cells meet across a common face. In particular $(abc)^5d$ commutes with $\alpha$, which is a rotation with the same axis. The image $h$ of $(abc)^5d$ in $L_2(q)$ therefore acts in the same way as a screw transformation on ${\cal I}/K$, commuting with $g=\overline\alpha$ and cyclically permuting the cells around a $\langle g\rangle$-invariant bracelet. It follows that $n_3$ is the least integer $n\geq 1$ such that $h^n\in\langle g\rangle$, that is, the order of the image of $h$ in $C(g)/\langle g\rangle$; here $C(g)$ denotes the centraliser of $g$ in $L_2(q)$, a cyclic group of order $(q\mp 1)/2$ or $15$ if $q\equiv\pm 1$ mod~$(6)$ or $q=2^4$, but isomorphic to $C_3\times C_3$ if $q=3^2$. The unique integer $t_3=0$ or $
\pm 1$ such that $h^{n_3}=g^{t_3}$ indicates the angle $2\pi t_3/3$ of twisting around the axis required to close up the sequence of $n_3$ cells to form the bracelet.

Now $(abc)^5d=(\alpha\beta^{-1}\alpha^{-1}\beta)^2\alpha\gamma$, so if $q$ is odd then $(3.1)$ and $(3.3)$ give
\[h=\frac{1}{2}\left(\begin{array}{cc} t(D+e) & ft+1-t \\ ft-1+t & t(D-e) \end{array}\right)\]
with trace $tD$, where $D=\sqrt{5+4t}$,
% Check this: we had $D=f^2(4t+5)$ in \S 3.2
while if $q=2^4$ then $(3.2)$ and $(3.3)$ give
\[h=\left(\begin{array}{cc} xt +t & 0 \\ t+1 & xt \end{array}\right).\]

In principle this enables one to compute the values of $n_3$ and $t_3$ in any specific case. Assume that $q\neq 3^2$, so $C(g)$ is cyclic. If the order $|h|$ of $h$ in $L_2(q)$ is coprime to $3$ then $n_3=|h|$ and $t_3=0$, whereas if $|h|$ is divisible by $3$ then $n_3=|h|/3$ and $t_3=\pm 1$. For instance, if $q=2^4$ then $|h|=5$, so $n_3=5$ and $t_3=0$; thus each $g$ leaves invariant a single bracelet of five cells, with no twisting. Computations in the cases where $q=p\equiv \pm 1$ mod~$(5)$ for primes $p\leq 151$ (see Table~1 in \S 11) show that $|h|=|C(g)|$ in most (though not all) of these cases, so that $g$ has one bracelet of $n_3=(q\mp 1)/6$ cells as $q\equiv \pm 1$ mod~$(6)$. However there are several exceptions, the first being $q=59$, where $|C(g)|=30$: by Example~4.3 there are two chiral pairs of manifolds $\cal M$ in this case (see also Example~10.3); one pair has $|h|=30$, so each $g$ has a single bracelet of $(q+1)/6=10$ cells, but the other pair has $|h|=15$, so each $g$ has two bracelets of five cells. There is a similar phenomenon for $q=71$, where there are also two chiral pairs: one has $|h|=36$, giving one bracelet of twelve cells, and the other has $h=12$, giving three bracelets of four cells.

\subsection{Elements of order $5$}

Now suppose that $g$ has order $5$. Similar arguments show that if $q\equiv\pm 1$ mod~$(5)$ then $g$ has $(q\mp 1)/10n_5$ axes of rotation, for some integer $n_5$, each enclosed by a `necklace' of $n_5$ $\langle g\rangle$-invariant cells and $n_5$ fixed vertices, where adjacent cells meet. The situation is similar if $q=2^4$, with $g$ having $3/n_5$ necklaces of $n_5$ cells and vertices. If $q=5^2$ then the elements of order $5$ conjugate to $\overline\beta$ have $10/n_5$ such necklaces, but those in the other class act on $\cal M$ without fixed points. As in the case of elements of order $3$, the axes for a given $g$ form an orbit under $N(\langle g\rangle)$, and those for all $g$ form an orbit under $L_2(q)$; in this case, however, the latter orbit is self-dual, that is, invariant under $\Omega^+$.

In order to evaluate $n_5$ we let $h\in L_2(q)$ be the image of the element $(abc)^5(bcd)^5=(\alpha\beta^{-1}\alpha^{-1}\beta)^2\alpha\beta^{-1}\gamma(\beta\gamma^{-1}\beta^{-1}\gamma)^2$ of $\Delta$, which translates a cell of $\cal I$ to a neighbour meeting it across a common vertex. Thus $h$ commutes with $g=\overline\beta$ and cyclically permutes the cells around a $\langle g\rangle$-invariant necklace, so $n_5$ is the order of the image of $h$ in $C(g)/\langle g\rangle$. In this case $C(g)$ is a cyclic group of order $(q\mp 1)/2$ or $15$ if $q\equiv\pm 1$ mod~$(10)$ or $q=2^4$, but isomorphic to $C_5\times C_5$ if $q=5^2$.

Since $(abc)^5(bcd)^5=(\alpha\beta^{-1}\alpha^{-1}\beta)^2\alpha\beta^{-1}\gamma(\beta\gamma^{-1}\beta^{-1}\gamma)^2$, $(3.1)$ and $(3.3)$ give
\[h=\frac{t}{2(2+t)^2}\left(\begin {array}{cc} -(2+t+fD(3+2t)) & D(3+2t)(1+e)\\ D(e+2t)(e-1) & (fD(3+2t)-(2+t)) \end{array}\right)\]
with trace $-t^3$, if $q$ is odd, while if $q=2^4$ then $(3.2)$ and $(3.3)$ give
\[h=\left(\begin{array}{cc} x+t & t \\ 1 & x+t+1 \end{array}\right).\]

As in the preceding section, one can use this in specific cases to compute $n_5$ and the twisting parameter $t_5=0, \pm 1$ or $\pm 2$, defined by $h^{n_5}=g^{t_5}$, where $g=\overline\beta$. If $q\neq 3^2$ then $C(g)$ is cyclic, so if $|h|$ is coprime to $5$ then $n_5=|h|$ and $t_5=0$, whereas if $|h|$ is divisible by $5$ then $n_5=|h|/5$ and $t_5\neq 0$.  For instance, if $q=2^4$ then $|h|=3$, so $n_5=3$ and $t_5=0$; thus each $g$ leaves invariant a single necklace of three cells, with no twisting. Again, computations for small primes $q=p\equiv \pm 1$ mod~$(5)$ show that $|h|=|C(g)|$ in many cases, so that $g$ has one necklace of $n_5=(q\mp 1)/10$ cells as $q\equiv \pm 1$ mod~$(10)$. The first exception is again $q=59$: here $|C(g)|=30$, whereas the two chiral pairs have $|h|=10$ and $15$, giving three necklaces of two cells or two necklaces of three cells respectively.

\subsection{Elements of order $2$}

Any element $g$ of order $2$ has axes of rotation whose passage through $\cal M$ is described by $n_2$ repetitions of the sequence vertex, face, edge, cell, edge, face. The number of such axes for each $g$ is $(q\mp 1)/n_2$ if $q\equiv \pm 1$ mod~$(4)$, or $4/n_2$ if $q=2^4$. In principle one could calculate $n_2$ and the twisting parameter $t_2=0$ or $1$ in any specific case by imitating the method used for elements of order $3$ or $5$. However the details are even more laborious, so we omit them.

\section{Subgroups acting freely on ${\cal M}={\cal H}^3/K$}

If $K$ is a normal subgroup of $\Delta$ with quotient $L_2(q)$ then there is an induced action of $L_2(q)$ as a group of orientation-preserving isometries of the manifold ${\cal M}={\cal H}^3/K$. If $S$ is any subgroup of $L_2(q)$ then ${\cal M}/S$ is an orbifold, which is a manifold if and only if $S$ acts freely (without fixed points) on $\cal M$. Now an element fixing a point in $\cal M$ must stabilise a vertex, edge, face or cell of ${\cal I}/K$; the stabilisers of these in $L_2(q)$ have orders $60, 6, 6$ and $60$ respectively, so if $|S|$ is coprime to $30$ then $S$ acts freely. The converse is also true if $q\neq 3^2$ or $5^2$, since in such cases it  follows from \S 8 that each element of order $2$, $3$ or $5$ has a fixed point. When $S$ acts freely the manifold ${\cal M}/S$ is tessellated by icosahedra, and the number of them is
\[\frac{|L_2(q):S|}{60}=\left\{
\begin{array}{ll}
q(q^2-1)/120|S| & \mbox{if $q$ is odd,}\\
&\\
68/|S| & \mbox{if $q=16$.}
\end{array}\right.
\]

In particular, if $|L_2(q):S|=60$, or equivalently if $S$ complements an icosahedral subgroup in $L_2(q)$, then ${\cal M}/S$ consists of a single icosahedron with its faces identified. We will show that such a subgroup exists in just three cases, namely $q=11$, $29$ and $59$. In~\cite{Eve}, Everitt has classified (up to isometry) the manifolds formed by identifying faces of a platonic solid; entries 9--14 in his Table~2, based on work of Richardson and Rubinstein~\cite{RR}, deal with the case of an icosahedron with dihedral angles $2\pi/3$, so ${\cal M}/S$ must correspond to one of these entries. We will show how to determine this entry by considering the homology of ${\cal M}/S$ and by finding the face-identifications.

If $S'$ denotes the commutator subgroup of $S$ then there is a regular abelian covering ${\cal M}/S'\to{\cal M}/S$, induced by the action of the abelianisation $S^{\rm ab}=S/S'$ as a group of isometries of ${\cal M}/S'$. It follows that $S^{\rm ab}$ is a quotient of the first integer homology group $H_1({\cal M}/S,{\bf Z})$ of ${\cal M}/S$. In~\cite{Eve} Everitt gives the homology groups of the manifolds he has classified, and in the cases $q=11$ and $29$ the fact that $S^{\rm ab}$ is a quotient is sufficient to determine the relevant entry in his Table~2. In the case $q=59$, however, there are two chiral pairs of manifolds ${\cal M}/S$, corresponding to two entries with isomorphic homology groups; to determine which pair corresponds to which entry we therefore need to determine the face-identifications, which are given in~\cite{Eve}. (This method, used in~\cite{JM} in the case $q=11$, can also be used to confirm the correspondences given by considering homology groups in the cases $q=11$ and $29$.)

To determine the face-identifications for ${\cal M}/S$, it is sufficient to consider the images under $S$ of the faces of a single cell in $\cal M$, or equivalently to see how $\tilde S=\theta^{-1}(S)$ transforms the faces of a cell of ${\cal I}=\{3,5,3\}$. The flags of $\cal I$ can be identified with the elements of $\Gamma$, permuted regularly by $\Gamma$ acting on the right as a group of monodromy permutations: each generator $a, b, c$ or $d$ of $\Gamma$ changes the vertex, edge, face or cell in a flag, so that each cell corresponds to a coset $g\Gamma_0$ of $\Gamma_0=\langle a, b, c\rangle$ in $\Gamma$. (This action of $\Gamma$ commutes with its action on the left as a group of isometries of $\cal I$.) It is convenient to use the cell whose flags correspond to the elements of $\Gamma_0$, or equivalently its image in ${\cal I}/K$. The identification of any face is determined by that of any of its flags, so in each face one can choose an `odd' flag, corresponding to some $g \in \Gamma_0 \setminus \Gamma_0^+$. The incident flag in the neighbouring cell then corresponds to the element $gd$ of $\Gamma^+$, whose image $h \in L_2(q)$ can be computed by writing $gd$ as a word in $\alpha, \beta, \gamma$ and then applying the appropriate epimorphism $\theta : \Gamma^+ \to L_2(q)$. The subgroup $S$ of $L_2(q)$ acts on flags on the left, by isometries of $\cal M$, so to find the flag of ${\cal M}/S$ corresponding to $gd$ we need the factorisation $h = st$ where $s \in S$ and $t\in I=\theta(\Gamma_0^+)<L_2(q)$: the element $s$ tells us which cell the flag $gd$ of ${\cal I}/K$ belongs to, and the element $t$ tells us which flag of that cell is identified with $gd$ under the action of $S$. This factorisation exists and is unique since $L_2(q)=SI$ with $S\cap I=1$, and one can find $s$ and $t$ by representing $I$ as a regular permutation group on the $60$ cosets of $S$ in $L_2(q)$:  given $h$ (now expressed as a permutation of these cosets) one can find the unique $t\in I$ which sends the coset $S$ to $Sh$, and then one computes $s=ht^{-1}$. Expressing $t$ as a word in the images of $\alpha=ab$ and $\beta=bc$ indicates which flag $S$ identifies with $gd$, and hence gives the identification of the chosen face of ${\cal M}/S$. Doing this for all 20 faces gives the complete identification.

\smallskip

\noindent{\bf Example 10.1.} Theorems~1.1(c), 1.2 and 1.3(c) show that if we take $q=11$, so that $K=K_0$, we obtain a manifold ${\cal M}={\cal H}^3/K_0$ with isometry group $\Omega/K\cong PGL_2(11)\times C_2$; this is the manifold ${\cal M}_0$ described in Theorem~A of~\cite{JM}. The subgroup $\Gamma/K\cong L_2(11)\times C_2$ preserves a tessellation of ${\cal M}$ by $|L_2(11)|/60=11$ icosahedra, and the quotient tessellation of the orbifold ${\cal M}/C_2$ is Coxeter's $11$-cell~\cite{Cox}, an abstract regular polytope with symmetry group $L_2(11)$. A Sylow $11$-subgroup $S\cong C_{11}$ of $L_2(11)$ acts freely on ${\cal M}$, permuting the eleven icosahedra regularly, so ${\cal M}/S$ is a manifold consisting of a single icosahedron with suitable identifications of its faces; as shown in~\cite{JM} it is a Fibonacci manifold, with the Fibonacci group $F(2,10)$ as its fundamental group. The normaliser of $S$ in $L_2(11)$ is a Borel subgroup $B$, with $B/S\cong C_5$, giving a $5$-fold cyclic covering ${\cal M}/S\to {\cal M}/B$; it follows from results of Helling, Kim and Mennicke~\cite{HKM} and of Hilden, Lozano and Montesinos~\cite{HLM} that this is a covering of the $3$-sphere $S^3$, branched over a figure-eight knot. It is also shown by Vesnin and Mednykh in~\cite{VM} that ${\cal M}/S$ is a double covering of $S^3$, branched over the closed $3$-braid $(\sigma_1\sigma_2^{-1})^5$.
% shown in Figure~1.

%\begin{center}
%\begin{figure}[ht]
% \begin{minipage}[b]{\linewidth}
%  \centering \includegraphics[width=1.5in, totalheight=1.5in]{knot1}
%  \caption{The closed braid $(\sigma_1\sigma_2^{-1})^5$} \label{fig1}
% \end{minipage}
%\end{figure}
%\end{center}
%%\centerline {Figure 1: closed braid (with $5$-fold symmetry)}

Since ${\cal M}/S$ has an $11$-fold regular covering ${\cal M}$, its first integer homology group $H_1({\cal M}/S,{\bf Z})$ has ${\bf Z}_{11}$ as a quotient. In Everitt's list of manifolds formed by identifying faces of an icosahedron with dihedral angles $2\pi/3$ (entries 9--14 in Table~2 of~\cite{Eve}), there is only one example where the homology group has a quotient ${\bf Z}_{11}$, namely entry~9 with homology group ${\bf Z}_{11}\oplus{\bf Z}_{11}$. It follows that ${\cal M}/S$ is this manifold, and ${\cal M}$ is an $11$-fold regular covering of it.

\medskip

\noindent{\bf Example 10.2.} Theorems~1.1(e) and 1.3(e) show that there are two normal subgroups $K$ of $\Delta$, conjugate in $\Gamma$, with $\Delta/K\cong L_2(29)$. These correspond to $t=5$, giving $\sqrt{4t+5}=\pm5$ and hence $s=5^2$ or $4^2$, so by \S 4 the corresponding chiral pair of manifolds ${\cal M}={\cal H}^3/K$ have isometry group $\Omega^+/K\cong L_2(29)\times C_2$. In each case, $L_2(29)$ preserves a tessellation by $203 = 7.29$ icosahedra (see the comment at the end of Section 5 of~\cite{JM}). Now $L_2(29)$ has a conjugacy class of subgroups $S$ of order $203$ which act freely on $\cal M$ and permute the icosahedra regularly: each $S$ is the unique subgroup of index 2 in a Borel subgroup, i.e.~in the normaliser of a Sylow $29$-subgroup $S'\cong C_{29}$ of $L_2(29)$. It follows that the manifold ${\cal M}/S$ consists of a single icosahedron with faces identified; it has a Klein four-group $C_2\times C_2$ of isometries, with one factor $C_2$ (induced by the Borel subgroup) preserving the icosahedron, and the other factor (induced by the direct factor $C_2$ in the isometry group) sending the icosahedron to its dual, which is also an icosahedron. This manifold has a $7$-sheeted regular covering ${\cal M}/S'$, and there is only one example in Everitt's list where the homology group has ${\bf Z}_7$ as a quotient, namely entry~12 with homology group ${\bf Z}_5\oplus{\bf Z}_7$. It follows that the chiral pair of manifolds ${\cal M}/S$ described here correspond to this entry, and the manifolds $\cal M$ are $203$-sheeted regular coverings of it.

We now show that, as in Example~10.1, the icosahedral manifold ${\cal M}/S$ is a double covering of $S^3$ branched over a knot.
%\begin{center}
%\begin{figure}[ht]
% \begin{minipage}[b]{\linewidth}
%  \centering \includegraphics[width=1.5in, totalheight=1.5in]{knot2}
%  \caption{The knot $\cal K$} \label{fig2}
% \end{minipage}
%\end{figure}
%\end{center}
%%\centerline {Figure 2: the knot $\cal K$}
First of all, Weeks's program SnapPea~\cite{Wee} shows that the double covering $\tilde{\cal I}$ of $S^3$ branched over
% the knot $\cal K$ in Figure~2
a certain knot $\cal K$ is a closed hyperbolic manifold, and that its fundamental group $\pi_1(\tilde{\cal I})$ has generators $U, V$ with defining relations
\begin{eqnarray*}
&&UV^{-1}U^{-2}V^{-1}UVU^{-1}VU^2VU^{-1}V^3U^{-1}VU^2VU^{-1}VUV^{-1}U^{-2}V^{-1}\\
&&\qquad=UVU^2VU^{-1}VU^2VU^{-1}V^2U^{-1}VU^2VU^{-1}VU^2V=1.
\end{eqnarray*}
(This presentation can also be obtained by hand from the Wirtinger presentation of the knot group $\pi_1(S^3\setminus{\cal K})$ corresponding to $\cal K$.) The fundamental group of ${\cal M}/S$ is the inverse image $\tilde S$ of $S$ in $\Delta$, and calculations with GAP show that $\tilde S$ has generators $u=bcda=bacd=\alpha^{-1}\gamma$ and $v=cbacdcbcba=\beta^{-1}\alpha\beta\gamma^{-1}\beta\alpha^{-1}$ with defining relations
\begin{eqnarray*}
&&uv^{-1}u^{-2}v^{-1}uvu^{-1}vu^2vu^{-1}v^3u^{-1}vu^2vu^{-1}vuv^{-1}u^{-2}v^{-1}\\
&&\qquad=uvu^2vu^{-1}vu^2vu^{-1}v^2u^{-1}vu^2vu^{-1}vu^2v=1.
\end{eqnarray*}
There is an obvious isomorphism $U\mapsto u, V\mapsto v$ between $\pi_1(\tilde{\cal I})$ and $\tilde S$. Since both $\tilde{\cal I}$ and ${\cal M}/S$ are closed hyperbolic manifolds, the Mostow Rigidity Theorem implies that they are isometric, so ${\cal M}/S$ is a double covering of $S^3$ branched over the knot $\cal K$.

\medskip

\noindent{\bf Example 10.3.} It follows from Example~4.3 that there are two chiral pairs of manifolds ${\cal M}={\cal H}^3/K$ corresponding to normal subgroups $K$ of $\Delta$ with $\Delta/K\cong L_2(59)$. For one pair the isometry group $\Omega^+/K$ is $PGL_2(59)$, for the other it is $L_2(59)\times C_2$. In each case the subgroup $L_2(59)$ preserves a tessellation by $1711=29.59$ icosahedra. A Borel subgroup $S$ of $L_2(59)$ (i.e.~the normaliser of a Sylow $59$-subgroup $S'\cong C_{59}$) acts freely on $\cal M$ and permutes the icosahedra regularly, so again the manifold ${\cal M}/S$ consists of a single icosahedron. The $29$-sheeted regular covering ${\cal M}/S'$ shows that $H_1({\cal M}/S,{\bf Z)}$ has ${\bf Z}_{29}$ as a quotient, so the two chiral pairs ${\cal M}/S$ correspond to entries~13 and 14 in Everitt's Table~2, both with homology group ${\bf Z}_{29}$. The manifolds $\cal M$ are $1711$-sheeted regular coverings of these.

%To determine which manifolds ${\cal M}/S$ correspond to entry~13 and which to entry~14, we need to find the face-identifications, using the general method described earlier, and to compare them with those given in~\cite{Eve}. We can choose the Borel subgroup $S$ to be the stabiliser of $0$ in the natural action of $L_2(59)$ by M\"obius transformations on the projective line $PG_1(59)$. If $h=st$ as before then $h^{-1}(0) = (st)^{-1}(0) = t^{-1}s^{-1}(0) = t^{-1}(0)$, where we compose M\"obius transformations from right to left, so that they are compatible with matrix multiplication in $SL_2(59)$. We therefore need to find the unique $t \in I$ such that $t^{-1}(0) = h^{-1}(0)$. This step is tedious: there are 60 elements $t\in I$ to check, this will need to be done for each of the 20 faces of the cell, and there are four manifolds $\cal M$ to consider. However, it is simple to program this process, and the resulting face-identifications indicate that the chiral pair with isometry group $PGL_2(59)$ correspond to Everitt's entry~?, while the pair with isometry group $L_2(59)\times C_2$ correspond to entry~??. We hope to investigate further geometric and topological properties of these manifolds ${\cal M}/S$, and of various other manifolds described here, in a future paper.

\medskip

Since $L_2(q)$ has no proper subgroup of index less than $q+1$ for $q>11$, $L_2(59)$ is the largest example of a group $L_2(q)$ in which an icosahedral subgroup has a complement. By inspection of the cases $q\leq 59$, using Dickson's classification~\cite[Ch.XII]{Dic} of the subgroups of $L_2(q)$, the only examples of this phenomenon (with a non-identity complement) are those described above with $q=11$, $29$ or $59$, so these are the only cases where $\cal M$ is a regular covering of a manifold consisting of a single icosahedron. There are, however, other examples where a quotient manifold ${\cal M}/S$ consists of a small number of icosahedra: for instance, taking $S$ to be a Sylow $19$-subgroup in $L_2(19)$ or a Sylow $17$-subgroup in $L_2(2^4)$ yields tessellations with three or four cells respectively.

\section{Tables of results}

By using GAP to check the relevant conditions on finite fields, we obtain the following tables of results. Successive columns indicate the prime $p$, the possible values of $\sqrt 5$ in $F_p$ (omitted when $p\equiv\pm 2$ mod~$(5)$) and of the trace $t_i=(-1\pm\sqrt 5)/2$, the resulting discriminant $D_i$ and the corresponding number of normal subgroups $K$. The next three columns give the quotients $\Gamma^+/K$, $\Gamma/K$ (with symbols  $\times$, $\Sigma$, $\Sigma^+$ and $-$ indicating that the quotient is $L_2(q)\times C_2$, $P\Sigma L_2(q)$, $P\Sigma L_2(q)^+$ or that $K$ is not normal in $\Gamma$), and $\Omega^+/K$ (with $\times$ and $\bullet$ indicating $L_2(q)\times C_2$ or $PGL_2(q)$). The final column indicates the corresponding case in Theorems~1.1 and 1.3.

\newpage

\begin{table}[!ht]
 \clearpage
 \centering
 \begin{tabular}{|c|c|l|l|c|c|c|c|c|}
 \hline
 $p$ & $\sqrt{5}$ &
 $t_1$, $t_2$ &
 $D_1$, $D_2$ &
 $\#$ & $\Gamma^+/K$ & $\Gamma/K$ & $\Omega^+/K$ & Case\\
 \hline
 \rowcolor[gray]{1} 5 & $0$ & $\pm 2$ & $\mp 2$ & $1$ & $L_2(5^5)$ & $\Sigma$ & $\times$ & (b)\\
 \rowcolor[gray]{0.9} 11 & $\pm 4$ &
  $\left\{\begin{array}{l} -4 \\ 3 \end{array}\right.$ &
  $\begin{array}{l} 0 \\ -5 \end{array}$ &
  $\begin{array}{l} 1 \\ 1 \end{array}$ &
  $\begin{array}{l} L_2(11) \\ L_2(11^2) \end{array}$ &
  $\begin{array}{l} \times \\ \Sigma \end{array}$ &
  $\begin{array}{l} \bullet \\ \times \end{array}$ & (c)\\
 \rowcolor[gray]{1} 19 & $\pm 10$ &
  $\left\{\begin{array}{l} 4 \\ -5 \end{array}\right.$ &
  $\begin{array}{l} 2 \\ 4=2^2 \end{array}$ &
  $\begin{array}{l} 1 \\ 2 \end{array}$ &
  $\begin{array}{l} L_2(19^2) \\ L_2(19) \end{array}$ &
  $\begin{array}{l} \Sigma \\ - \end{array}$ &
  $\begin{array}{l} \times \\ \bullet \end{array}$ & (e)\\
 \rowcolor[gray]{0.9} 29 & $\pm 11$ &
  $\left\{\begin{array}{l} 5 \\ -6 \end{array}\right.$ &
  $\begin{array}{l} -4=5^2 \\ 10 \end{array}$ &
  $\begin{array}{l} 2 \\ 1 \end{array}$ &
  $\begin{array}{l} L_2(29) \\ L_2(29^2) \end{array}$ &
  $\begin{array}{l} - \\ \Sigma \end{array}$ &
  $\begin{array}{l} \times \\ \times \end{array}$ & (e)\\
 \rowcolor[gray]{1} 31 & $\pm 6$ &
  $\left\{\begin{array}{l} -13 \\ 12 \end{array}\right.$ &
  $\begin{array}{l} -16 \\ -9 \end{array}$ &
  $\begin{array}{l} 1 \\ 1 \end{array}$ &
  $\begin{array}{l} L_2(31^2) \\ L_2(31^2) \end{array}$ &
  $\begin{array}{l} \Sigma \\ \Sigma \end{array}$ &
  $\begin{array}{l} \times \\ \times \end{array}$ & (d)\\
 \rowcolor[gray]{0.9} 41 & $\pm 13$ &
  $\left\{\begin{array}{l} 6 \\ -7 \end{array}\right.$ &
  $\begin{array}{l} -12 \\ 18=10^2 \end{array}$ &
  $\begin{array}{l} 1 \\ 2 \end{array}$ &
  $\begin{array}{l} L_2(41^2) \\ L_2(41) \end{array}$ &
  $\begin{array}{l} \Sigma \\ - \end{array}$ &
  $\begin{array}{l} \times \\ \bullet \end{array}$ & (e)\\
 \rowcolor[gray]{1} 59 & $\pm 8$ &
  $\left\{\begin{array}{l} -26 \\ 25 \end{array}\right.$ &
  $\begin{array}{l} 19 \\ -13 \end{array}$ &
  $\begin{array}{l} 2 \\ 2 \end{array}$ &
  $\begin{array}{l} L_2(59) \\ L_2(59) \end{array}$ &
  $\begin{array}{l} - \\ - \end{array}$ &
  $\begin{array}{l} \bullet \\ \times \end{array}$ & (d)\\
 \rowcolor[gray]{0.9} 61 & $\pm 26$ &
  $\left\{\begin{array}{l} -18 \\ 17 \end{array}\right.$ &
  $\begin{array}{l} -6 \\ 12=16^2 \end{array}$ &
  $\begin{array}{l} 1 \\ 2 \end{array}$ &
  $\begin{array}{l} L_2(61^2) \\ L_2(61) \end{array}$ &
  $\begin{array}{l} \Sigma \\ - \end{array}$ &
  $\begin{array}{l} \times \\ \bullet \end{array}$ & (e)\\
 \rowcolor[gray]{1} 71 & $\pm 17$ &
  $\left\{\begin{array}{l} 8 \\ -9 \end{array}\right.$ &
  $\begin{array}{l} -34=26^2 \\ -31=18^2 \end{array}$ &
  $\begin{array}{l} 2 \\ 2 \end{array}$ &
  $\begin{array}{l} L_2(71) \\ L_2(71) \end{array}$ &
  $\begin{array}{l} - \\ - \end{array}$ &
  $\begin{array}{l} \times \\ \bullet \end{array}$ & (d)\\
 \rowcolor[gray]{0.9} 79 & $\pm 20$ &
  $\left\{\begin{array}{l} -30 \\ 29 \end{array}\right.$ &
  $\begin{array}{l} -36 \\ -37=11^2 \end{array}$ &
  $\begin{array}{l} 1 \\ 2 \end{array}$ &
  $\begin{array}{l} L_2(79^2) \\ L_2(79) \end{array}$ &
  $\begin{array}{l} \Sigma \\ - \end{array}$ &
  $\begin{array}{l} \times \\ \bullet \end{array}$ & (e)\\
 \rowcolor[gray]{1} 89 & $\pm 19$ &
  $\left\{\begin{array}{l} 9 \\ -10 \end{array}\right.$ &
  $\begin{array}{l} 41 \\ -35 \end{array}$ &
  $\begin{array}{l} 1 \\ 1 \end{array}$ &
  $\begin{array}{l} L_2(89^2) \\ L_2(89^2) \end{array}$ &
  $\begin{array}{l} \Sigma \\ \Sigma \end{array}$ &
  $\begin{array}{l} \times \\ \times \end{array}$ & (d)\\
 \rowcolor[gray]{0.9} 101 & $\pm 45$ &
  $\left\{\begin{array}{l} 22 \\ -23 \end{array}\right.$ &
  $\begin{array}{l} -8 \\ 14=32^2 \end{array}$ &
  $\begin{array}{l} 1 \\ 2 \end{array}$ &
  $\begin{array}{l} L_2(101^2) \\ L_2(101) \end{array}$ &
  $\begin{array}{l} \Sigma \\ - \end{array}$ &
  $\begin{array}{l} \times \\ \times \end{array}$ & (e)\\
 \rowcolor[gray]{1} 109 & $\pm 21$ &
  $\left\{\begin{array}{l} 10 \\ -11 \end{array}\right.$ &
  $\begin{array}{l} 45=46^2 \\ -39 \end{array}$ &
  $\begin{array}{l} 2 \\ 1 \end{array}$ &
  $\begin{array}{l} L_2(109) \\ L_2(109^2) \end{array}$ &
  $\begin{array}{l} - \\ \Sigma \end{array}$ &
  $\begin{array}{l} \bullet \\ \times \end{array}$ & (e)\\
 \rowcolor[gray]{0.9} 131 & $\pm 23$ &
  $\left\{\begin{array}{l} 11 \\ -12 \end{array}\right.$ &
  $\begin{array}{l} -43 \\ 49=7^2 \end{array}$ &
  $\begin{array}{l} 1 \\ 2 \end{array}$ &
  $\begin{array}{l} L_2(131^2) \\ L_2(131) \end{array}$ &
  $\begin{array}{l} \Sigma \\ - \end{array}$ &
  $\begin{array}{l} \times \\ \times \end{array}$ & (e)\\
 \rowcolor[gray]{1} 139 & $\pm 12$ &
  $\left\{\begin{array}{l} 75 \\ 63 \end{array}\right.$ &
  $\begin{array}{l} -21=37^2 \\ 27 \end{array}$ &
  $\begin{array}{l} 2 \\ 1 \end{array}$ &
  $\begin{array}{l} L_2(139) \\ L_2(139^2) \end{array}$ &
  $\begin{array}{l} - \\ \Sigma \end{array}$ &
  $\begin{array}{l} \times \\ \times \end{array}$ & (e)\\
 \rowcolor[gray]{0.9} 149 & $\pm 68$ &
  $\left\{\begin{array}{l} 108 \\ 40 \end{array}\right.$ &
  $\begin{array}{l} -10 \\ 16=4^2 \end{array}$ &
  $\begin{array}{l} 1 \\ 2 \end{array}$ &
  $\begin{array}{l} L_2(14^2) \\ L_2(149) \end{array}$ &
  $\begin{array}{l} \Sigma \\ - \end{array}$ &
  $\begin{array}{l} \times \\ \bullet \end{array}$ & (e)\\
 \rowcolor[gray]{1} 151 & $\pm 55$ &
  $\left\{\begin{array}{l} 27 \\ -28 \end{array}\right.$ &
  $\begin{array}{l} 44=73^2 \\ -38 \end{array}$ &
  $\begin{array}{l} 2 \\ 1 \end{array}$ &
  $\begin{array}{l} L_2(151) \\ L_2(151^2) \end{array}$ &
  $\begin{array}{l} - \\ \Sigma \end{array}$ &
  $\begin{array}{l} \times \\ \times \end{array}$ & (e)\\
 \rowcolor[gray]{0.9} 179 & $\pm 30$ &
  $\left\{\begin{array}{l} 104 \\ 74 \end{array}\right.$ &
  $\begin{array}{l} -57 \\ 63 \end{array}$ &
  $\begin{array}{l} 1 \\ 1 \end{array}$ &
  $\begin{array}{l} L_2(179^2) \\ L_2(179^2) \end{array}$ &
  $\begin{array}{l} \Sigma \\ \Sigma \end{array}$ &
  $\begin{array}{l} \times \\ \times \end{array}$ & (d)\\
 \rowcolor[gray]{1} 181 & $\pm 27$ &
  $\left\{\begin{array}{l} 13 \\ 14 \end{array}\right.$ &
  $\begin{array}{l} 57 \\ -51 \end{array}$ &
  $\begin{array}{l} 1 \\ 1 \end{array}$ &
  $\begin{array}{l} L_2(181^2) \\ L_2(181^2) \end{array}$ &
  $\begin{array}{l} \Sigma \\ \Sigma \end{array}$ &
  $\begin{array}{l} \times \\ \times \end{array}$ & (d)\\
 \rowcolor[gray]{0.9} 191 & $\pm 14$ &
  $\left\{\begin{array}{l} 102 \\ 88 \end{array}\right.$ &
  $\begin{array}{l} -25 \\ 31 \end{array}$ &
  $\begin{array}{l} 1 \\ 1 \end{array}$ &
  $\begin{array}{l} L_2(191^2) \\ L_2(191^2) \end{array}$ &
  $\begin{array}{l} \Sigma \\ \Sigma \end{array}$ &
  $\begin{array}{l} \times \\ \times \end{array}$ & (d)\\
 \hline
 \end{tabular}
 \caption{The cases $p\equiv 0, \pm 1$ mod 5}
 \label{table_1a}
\end{table}

\newpage

\begin{table}[!ht]
 \clearpage
 \centering
 \begin{tabular}{|c|c|l|l|c|c|c|c|c|}
 \hline
 $p$ & $\sqrt{5}$ &
 $t_1$, $t_2$ &
 $D_1$, $D_2$ &
 $\#$ & $\Gamma^+/K$ & $\Gamma/K$ & $\Omega^+/K$ & Case\\
 \hline
 \rowcolor[gray]{1} 199 & $\pm 76$ &
  $\left\{\begin{array}{l} 137 \\ 61 \end{array}\right.$ &
  $\begin{array}{l} 50=100^2 \\ -44=73^2 \end{array}$ &
  $\begin{array}{l} 2 \\ 2 \end{array}$ &
  $\begin{array}{l} L_2(199) \\ L_2(199) \end{array}$ &
  $\begin{array}{l} - \\ - \end{array}$ &
  $\begin{array}{l} \times \\ \bullet \end{array}$ & (d)\\
 \rowcolor[gray]{0.9} 211 & $\pm 65$ &
  $\left\{\begin{array}{l} 32 \\ -33 \end{array}\right.$ &
  $\begin{array}{l} 84=57^2 \\ -78 \end{array}$ &
  $\begin{array}{l} 2 \\ 1 \end{array}$ &
  $\begin{array}{l} L_2(211) \\ L_2(211^2) \end{array}$ &
  $\begin{array}{l} - \\ \Sigma \end{array}$ &
  $\begin{array}{l} \times \\ \times \end{array}$ & (e)\\
 \rowcolor[gray]{1} 229 & $\pm 66$ &
  $\left\{\begin{array}{l} 147 \\ 81 \end{array}\right.$ &
  $\begin{array}{l} -94=89^2 \\ 100=10^2 \end{array}$ &
  $\begin{array}{l} 2 \\ 2 \end{array}$ &
  $\begin{array}{l} L_2(229) \\ L_2(229) \end{array}$ &
  $\begin{array}{l} - \\ - \end{array}$ &
  $\begin{array}{l} \times \\ \times \end{array}$ & (d)\\
 \rowcolor[gray]{0.9} 239 & $\pm 31$ &
  $\left\{\begin{array}{l} 15 \\ -16 \end{array}\right.$ &
  $\begin{array}{l} 65 \\ -59=53^2 \end{array}$ &
  $\begin{array}{l} 1 \\ 2 \end{array}$ &
  $\begin{array}{l} L_2(239^2) \\ L_2(239) \end{array}$ &
  $\begin{array}{l} \Sigma \\ - \end{array}$ &
  $\begin{array}{l} \times \\ \bullet \end{array}$ & (e)\\
 \rowcolor[gray]{1} 241 & $\pm 103$ &
  $\left\{\begin{array}{l} 51 \\ -52 \end{array}\right.$ &
  $\begin{array}{l} 38 \\ -32=89^2 \end{array}$ &
  $\begin{array}{l} 1 \\ 2 \end{array}$ &
  $\begin{array}{l} L_2(241^2) \\ L_2(241) \end{array}$ &
  $\begin{array}{l} \Sigma \\ - \end{array}$ &
  $\begin{array}{l} \times \\ \bullet \end{array}$ & (e)\\
 \rowcolor[gray]{0.9} 251 & $\pm 16$ &
  $\left\{\begin{array}{l} 133 \\ 117 \end{array}\right.$ &
  $\begin{array}{l} 35=81^2 \\ -29=67^2 \end{array}$ &
  $\begin{array}{l} 2 \\ 2 \end{array}$ &
  $\begin{array}{l} L_2(251) \\ L_2(251) \end{array}$ &
  $\begin{array}{l} - \\ - \end{array}$ &
  $\begin{array}{l} \bullet \\ \times \end{array}$ & (d)\\
 \hline
 \end{tabular}
 \caption{The case $p\equiv\pm 1$ mod 5, continued}
 \label{table_1b}
\end{table}

\begin{table}[!ht]
 \clearpage
 \centering
 \begin{tabular}{|c|l|l|c|c|c|c|c|}
 \hline
 $p$ &
 $t_1$, $t_2$ &
 $D_1$, $D_2$ &
 $\#$ & $\Gamma^+/K$ & $\Gamma/K$ & $\Omega^+/K$ & Case\\
 \hline
 \rowcolor[gray]{0.9} 2 &
  $\left\{\begin{array}{l} t_1 \\ t_2 \end{array}\right.$ &
  $\left.\begin{array}{l} \sqrt{5+4t_1} \\ \sqrt{5+4t_2} \end{array}\right\}$ &
  $1$ &
  $L_2(2^4)$ & $\Sigma^+$ & $\times$ & (a) \\
 \rowcolor[gray]{1} 3 &
  $\left\{\begin{array}{l} t_1 \\ t_2 \end{array}\right.$ &
  $\left.\begin{array}{l} 1+t_1 \\ 1+t_2 \end{array}\right\}$ &
  $\begin{array}{l} 1 \\ 1 \end{array}$ &
  $\begin{array}{l} L_2(3^2) \\ L_2(3^2) \end{array}$ &
  $\begin{array}{l} - \\ - \end{array}$ &
  $\begin{array}{l} \bullet \\ \bullet \end{array}$ & (f)\\
 \rowcolor[gray]{0.9} 7 &
  $\left\{\begin{array}{l} t_1 \\ t_2 \end{array}\right.$ &
  $\left.\begin{array}{l} \sqrt{5+4t_1} \\ \sqrt{5+4t_2} \end{array}\right\}$ &
  $1$ &
  $L_2(7^4)$ & $\Sigma^+$ & $\times$ & (g) \\
 \rowcolor[gray]{1} 13 &
  $\left\{\begin{array}{l} t_1 \\ t_2 \end{array}\right.$ &
  $\left.\begin{array}{l} \sqrt{5+4t_1} \\ \sqrt{5+4t_2} \end{array}\right\}$ &
  $1$ &
  $L_2(13^4)$ & $\Sigma^+$ & $\times$ & (g) \\
 \rowcolor[gray]{0.9} 17 &
  $\left\{\begin{array}{l} t_1 \\ t_2 \end{array}\right.$ &
  $\left.\begin{array}{l} \sqrt{5+4t_1} \\ \sqrt{5+4t_2} \end{array}\right\}$ &
  $1$ &
  $L_2(17^4)$ & $\Sigma^+$ & $\times$ & (g) \\
 \rowcolor[gray]{1} 23 &
  $\left\{\begin{array}{l} t_1 \\ t_2 \end{array}\right.$ &
  $\left.\begin{array}{l} 9+19t_1 \\ 9+19t_2 \end{array}\right\}$ &
  $\begin{array}{l} 1 \\ 1 \end{array}$ &
  $\begin{array}{l} L_2(23^2) \\ L_2(23^2) \end{array}$ &
  $\begin{array}{l} - \\ - \end{array}$ &
  $\begin{array}{l} \bullet \\ \bullet \end{array}$ & (f)\\
 \rowcolor[gray]{0.9} 37 &
  $\left\{\begin{array}{l} t_1 \\ t_2 \end{array}\right.$ &
  $\left.\begin{array}{l} 20+30t_1 \\ 17+7t_2 \end{array}\right\}$ &
  $\begin{array}{l} 1 \\ 1 \end{array}$ &
  $\begin{array}{l} L_2(37^2) \\ L_2(37^2) \end{array}$ &
  $\begin{array}{l} - \\ - \end{array}$ &
  $\begin{array}{l} \times \\ \times \end{array}$ & (f)\\
 \rowcolor[gray]{1} 43 &
  $\left\{\begin{array}{l} t_1 \\ t_2 \end{array}\right.$ &
  $\left.\begin{array}{l} \sqrt{5+4t_1} \\ \sqrt{5+4t_2} \end{array}\right\}$ &
  $1$ &
  $L_2(43^4)$ & $\Sigma^+$ & $\times$ & (g) \\
 \rowcolor[gray]{0.9} 47 &
  $\left\{\begin{array}{l} t_1 \\ t_2 \end{array}\right.$ &
  $\left.\begin{array}{l} 42+11t_1 \\ 42+11t_2 \end{array}\right\}$ &
  $\begin{array}{l} 1 \\ 1 \end{array}$ &
  $\begin{array}{l} L_2(47^2) \\ L_2(47^2) \end{array}$ &
  $\begin{array}{l} - \\ - \end{array}$ &
  $\begin{array}{l} \bullet \\ \bullet \end{array}$ & (f)\\
 \hline
 \end{tabular}
 \caption{The case $p\equiv\pm 2$ mod 5}
 \label{table_2a}
\end{table}

\newpage

\begin{table}[!ht]
 \clearpage
 \centering
 \begin{tabular}{|c|l|l|c|c|c|c|c|}
 \hline
 $p$ &
 $t_1$, $t_2$ &
 $D_1$, $D_2$ &
 $\#$ & $\Gamma^+/K$ & $\Gamma/K$ & $\Omega^+/K$ & Case\\
 \hline
 \rowcolor[gray]{1} 53 &
  $\left\{\begin{array}{l} t_1 \\ t_2 \end{array}\right.$ &
  $\begin{array}{l} 11 + 19t_1 \\ 11 + 19t_2 \end{array}$ &
  $\begin{array}{l} 1 \\ 1 \end{array}$ &
  $\begin{array}{l} L_2(53^2) \\ L_2(53^2) \end{array}$ &
  $\begin{array}{l} - \\ - \end{array}$ &
  $\begin{array}{l} \times \\ \times \end{array}$ & (f)\\
 \rowcolor[gray]{0.9} 67 &
  $\left\{\begin{array}{l} t_1 \\ t_2 \end{array}\right.$ &
  $\begin{array}{l} 12 + 53t_1 \\ 12 + 53t_2 \end{array}$ &
  $\begin{array}{l} 1 \\ 1 \end{array}$ &
  $\begin{array}{l} L_2(67^2) \\ L_2(67^2) \end{array}$ &
  $\begin{array}{l} - \\ - \end{array}$ &
  $\begin{array}{l} \bullet \\ \bullet \end{array}$ & (f)\\
 \rowcolor[gray]{1} 73 &
  $\left\{\begin{array}{l} t_1 \\ t_2 \end{array}\right.$ &
  $\left.\begin{array}{l} \sqrt{5+4t_1} \\ \sqrt{5+4t_2} \end{array}\right\}$ &
  $1$ &
  $L_2(73^4)$ & $\Sigma^+$ & $\times$ & (f) \\
 \rowcolor[gray]{0.9} 83 &
  $\left\{\begin{array}{l} t_1 \\ t_2 \end{array}\right.$ &
  $\left.\begin{array}{l} \sqrt{5+4t_1} \\ \sqrt{5+4t_2} \end{array}\right\}$ &
  $1$ &
  $L_2(83^4)$ & $\Sigma^+$ & $\times$ & (g) \\
 \rowcolor[gray]{1} 97 &
  $\left\{\begin{array}{l} t_1 \\ t_2 \end{array}\right.$ &
  $\begin{array}{l} 28 + 26t_1 \\ 28 + 26t_2 \end{array}$ &
  $\begin{array}{l} 1 \\ 1 \end{array}$ &
  $\begin{array}{l} L_2(97^2) \\ L_2(97^2) \end{array}$ &
  $\begin{array}{l} - \\ - \end{array}$ &
  $\begin{array}{l} \times \\ \times \end{array}$ & (f)\\
 \rowcolor[gray]{0.9} 103 &
  $\left\{\begin{array}{l} t_1 \\ t_2 \end{array}\right.$ &
  $\begin{array}{l} 33 + 96t_1 \\ 70 + 7t_2 \end{array}$ &
  $\begin{array}{l} 1 \\ 1 \end{array}$ &
  $\begin{array}{l} L_2(103^2) \\ L_2(103^2) \end{array}$ &
  $\begin{array}{l} - \\ - \end{array}$ &
  $\begin{array}{l} \bullet \\ \bullet \end{array}$ & (f)\\
 \rowcolor[gray]{1} 107 &
  $\left\{\begin{array}{l} t_1 \\ t_2 \end{array}\right.$ &
  $\left.\begin{array}{l} \sqrt{5+4t_1} \\ \sqrt{5+4t_2} \end{array}\right\}$ &
  $1$ &
  $L_2(107^4)$ & $\Sigma^+$ & $\times$ & (g) \\
 \rowcolor[gray]{0.9} 113 &
  $\left\{\begin{array}{l} t_1 \\ t_2 \end{array}\right.$ &
  $\begin{array}{l} 73 + 10t_1 \\ 73 + 10t_2 \end{array}$ &
  $\begin{array}{l} 1 \\ 1 \end{array}$ &
  $\begin{array}{l} L_2(113^2) \\ L_2(113^2) \end{array}$ &
  $\begin{array}{l} - \\ - \end{array}$ &
  $\begin{array}{l} \times \\ \times \end{array}$ & (f)\\
 \rowcolor[gray]{1} 127 &
  $\left\{\begin{array}{l} t_1 \\ t_2 \end{array}\right.$ &
  $\left.\begin{array}{l} \sqrt{5+4t_1} \\ \sqrt{5+4t_2} \end{array}\right\}$ &
  $1$ &
  $L_2(127^4)$ & $\Sigma^+$ & $\times$ & (g) \\
 \rowcolor[gray]{0.9} 137 &
  $\left\{\begin{array}{l} t_1 \\ t_2 \end{array}\right.$ &
  $\begin{array}{l} 30 + 129t_1 \\ 30 + 129t_2 \end{array}$ &
  $\begin{array}{l} 1 \\ 1 \end{array}$ &
  $\begin{array}{l} L_2(137^2) \\ L_2(137^2) \end{array}$ &
  $\begin{array}{l} - \\ - \end{array}$ &
  $\begin{array}{l} \times \\ \times \end{array}$ & (f)\\
 \rowcolor[gray]{1} 157 &
  $\left\{\begin{array}{l} t_1 \\ t_2 \end{array}\right.$ &
  $\begin{array}{l} 43 + 121t_1 \\ 114 + 36t_2 \end{array}$ &
  $\begin{array}{l} 1 \\ 1 \end{array}$ &
  $\begin{array}{l} L_2(157^2) \\ L_2(157^2) \end{array}$ &
  $\begin{array}{l} - \\ - \end{array}$ &
  $\begin{array}{l} \times \\ \times \end{array}$ & (f)\\
 \rowcolor[gray]{0.9} 163 &
  $\left\{\begin{array}{l} t_1 \\ t_2 \end{array}\right.$ &
  $\begin{array}{l} 12 + 137t_1 \\ 151 + 26t_2 \end{array}$ &
  $\begin{array}{l} 1 \\ 1 \end{array}$ &
  $\begin{array}{l} L_2(163^2) \\ L_2(163^2) \end{array}$ &
  $\begin{array}{l} - \\ - \end{array}$ &
  $\begin{array}{l} \bullet \\ \bullet \end{array}$ & (f)\\
 \rowcolor[gray]{1} 173 &
  $\left\{\begin{array}{l} t_1 \\ t_2 \end{array}\right.$ &
  $\left.\begin{array}{l} \sqrt{5+4t_1} \\ \sqrt{5+4t_2} \end{array}\right\}$ &
  $1$ &
  $L_2(173^4)$ & $\Sigma^+$ & $\times$ & (g) \\
 \rowcolor[gray]{0.9} 193 &
  $\left\{\begin{array}{l} t_1 \\ t_2 \end{array}\right.$ &
  $\left.\begin{array}{l} \sqrt{5+4t_1} \\ \sqrt{5+4t_2} \end{array}\right\}$ &
  $1$ &
  $L_2(193^4)$ & $\Sigma^+$ & $\times$ & (g) \\
 \rowcolor[gray]{1} 197 &
  $\left\{\begin{array}{l} t_1 \\ t_2 \end{array}\right.$ &
  $\left.\begin{array}{l} \sqrt{5+4t_1} \\ \sqrt{5+4t_2} \end{array}\right\}$ &
  $1$ &
  $L_2(197^4)$ & $\Sigma^+$ & $\times$ & (g) \\
 \hline
 \end{tabular}
 \caption{The case $p\equiv\pm 2$ mod 5, continued}
 \label{table_2b}
\end{table}

\newpage

\begin{table}[!ht]
 \clearpage
 \centering
 \begin{tabular}{|c|l|l|c|c|c|c|c|}
 \hline
 $p$ &
 $t_1$, $t_2$ &
 $D_1$, $D_2$ &
 $\#$ & $\Gamma^+/K$ & $\Gamma/K$ & $\Omega^+/K$ & Case\\
 \hline
 \rowcolor[gray]{0.9} 127 &
  $\left\{\begin{array}{l} t_1 \\ t_2 \end{array}\right.$ &
  $\left.\begin{array}{l} \sqrt{5+4t_1} \\ \sqrt{5+4t_2} \end{array}\right\}$ &
  $1$ &
  $L_2(127^4)$ & $\Sigma^+$ & $\times$ & (g) \\
 \rowcolor[gray]{1} 223 &
  $\left\{\begin{array}{l} t_1 \\ t_2 \end{array}\right.$ &
  $\left.\begin{array}{l} 118+199t_1 \\ 105+24t_2 \end{array}\right\}$ &
  $\begin{array}{l} 1 \\ 1 \end{array}$ &
  $\begin{array}{l} L_2(223^2) \\ L_2(223^2) \end{array}$ &
  $\begin{array}{l} - \\ - \end{array}$ &
  $\begin{array}{l} \bullet \\ \bullet \end{array}$ & (f)\\
 \rowcolor[gray]{0.9} 227 &
  $\left\{\begin{array}{l} t_1 \\ t_2 \end{array}\right.$ &
  $\left.\begin{array}{l} \sqrt{5+4t_1} \\ \sqrt{5+4t_2} \end{array}\right\}$ &
  $1$ &
  $L_2(227^4)$ & $\Sigma^+$ & $\times$ & (g) \\
 \rowcolor[gray]{1} 233 &
  $\left\{\begin{array}{l} t_1 \\ t_2 \end{array}\right.$ &
  $\left.\begin{array}{l} \sqrt{5+4t_1} \\ \sqrt{5+4t_2} \end{array}\right\}$ &
  $1$ &
  $L_2(233^4)$ & $\Sigma^+$ & $\times$ & (g) \\
 \hline
 \end{tabular}
 \caption{The case $p\equiv\pm 2$ mod 5, continued}
 \label{table_2c}
\end{table}

\end{document}